\documentclass[12pt,leqno]{article}
\usepackage{amssymb}
\usepackage[mathscr]{eucal}
\usepackage{amsmath,amssymb,latexsym,theorem,bbm}
\usepackage{amsmath,amssymb,latexsym,theorem}
\usepackage{color}
\usepackage{appendix}

\newcommand{\SC}{\scriptstyle}

\newcommand{\CC}{\mathbb{C}}
\newcommand{\DD}{\mathbb{D}}

\newcommand{\NN}{\mathbb{N}}

\newcommand{\RR}{\mathbb{R}}

\newcommand{\ZZ}{\mathbb{Z}}

\newcommand{\bA}{{\boldsymbol{A}}}

\newcommand{\ta}{\widetilde{a}}

\newcommand{\tbA}{\widetilde{\bA}}

\newcommand{\bB}{{\boldsymbol{B}}}

\newcommand{\tb}{\widetilde{b}}

\newcommand{\tbB}{\widetilde{\bB}}

\newcommand{\bc}{{\boldsymbol{c}}}
\newcommand{\bC}{{\boldsymbol{C}}}
\newcommand{\obC}{\overline{\boldsymbol{C}}}

\newcommand{\tbC}{\widetilde{\bC}}

\newcommand{\tc}{\widetilde{c}}

\newcommand{\be}{{\boldsymbol{e}}}
\newcommand{\tbe}{\widetilde{\be}}
\newcommand{\tf}{\widetilde{f}}
\newcommand{\Bf}{{\boldsymbol{f}}}

\newcommand{\bI}{{\boldsymbol{I}}}

\newcommand{\bJ}{{\boldsymbol{J}}}

\newcommand{\bM}{{\boldsymbol{M}}}

\newcommand{\oC}{\overline{C}}

\newcommand{\bP}{{\boldsymbol{P}}}

\newcommand{\bq}{{\boldsymbol{q}}}
\newcommand{\bQ}{{\boldsymbol{Q}}}

\newcommand{\bcQ}{{\boldsymbol{\cQ}}}

\newcommand{\bcR}{{\boldsymbol{\cR}}}
\newcommand{\bS}{{\boldsymbol{S}}}
\newcommand{\bu}{{\boldsymbol{u}}}

\newcommand{\bv}{{\boldsymbol{v}}}
\newcommand{\bV}{{\boldsymbol{V}}}

\newcommand{\bx}{{\boldsymbol{x}}}
\newcommand{\bX}{{\boldsymbol{X}}}
\newcommand{\by}{{\boldsymbol{y}}}

\newcommand{\bz}{{\boldsymbol{z}}}

\newcommand{\bU}{{\boldsymbol{U}}}

\newcommand{\Bbeta}{{\boldsymbol{\beta}}}
\newcommand{\tbeta}{\widetilde{\beta}}
\newcommand{\tBbeta}{\widetilde{\Bbeta}}

\newcommand{\ttBbeta}{\widetilde{\tBbeta}}
\newcommand{\bgamma}{{\boldsymbol{\gamma}}}

\newcommand{\tlambda}{\widetilde{\lambda}}

\newcommand{\bmu}{{\boldsymbol{\mu}}}

\newcommand{\bbeta}{{\boldsymbol{\beta}}}

\newcommand{\bPi}{{\boldsymbol{\Pi}}}

\newcommand{\bzero}{{\boldsymbol{0}}}
\newcommand{\bone}{{\boldsymbol{1}}}
\newcommand{\tbone}{\widetilde{\bone}}

\newcommand{\cA}{{\mathcal A}}
\newcommand{\cB}{{\mathcal B}}

\newcommand{\cD}{{\mathcal D}}

\newcommand{\cF}{{\mathcal F}}

\newcommand{\cM}{{\mathcal M}}
\newcommand{\bcM}{\boldsymbol{\cM}}

\newcommand{\cP}{{\mathcal P}}
\newcommand{\cQ}{{\mathcal Q}}
\newcommand{\cR}{{\mathcal R}}

\newcommand{\cU}{{\mathcal U}}

\newcommand{\bcU}{\boldsymbol{\cU}}

\newcommand{\cX}{{\mathcal X}}

\newcommand{\bcX}{\boldsymbol{\cX}}

\newcommand{\cY}{{\mathcal Y}}

\newcommand{\cW}{{\mathcal W}}
\newcommand{\bcW}{\boldsymbol{\cW}}
\newcommand{\bcY}{\boldsymbol{\cY}}

\newcommand{\cc}{\mathrm{c}}
\newcommand{\dd}{\mathrm{d}}
\newcommand{\ee}{\mathrm{e}}

\newcommand{\slu}{{\SC\mathrm{lu}}}

\newcommand{\INARp}{\textup{INAR($p$)}}

\newcommand{\EE}{\operatorname{\mathbb{E}}}

\newcommand{\PP}{\operatorname{\mathbb{P}}}
\newcommand{\OO}{\operatorname{O}}

\newcommand{\var}{\operatorname{Var}}

\renewcommand{\Re}{\operatorname{Re}}

\renewcommand{\mid}{\,|\,}
\newcommand{\bmid}{\,\big|\,}

\renewcommand{\leq}{\leqslant}
\renewcommand{\geq}{\geqslant}

\newcommand{\stoch}{\stackrel{\PP}{\longrightarrow}}
\newcommand{\distr}{\stackrel{\cD}{\longrightarrow}}

\newcommand{\lu}{\stackrel{\slu}{\longrightarrow}}

\newcommand{\bbone}{\mathbbm{1}}
\newcommand{\ns}{{\lfloor ns\rfloor}}
\newcommand{\nt}{{\lfloor nt\rfloor}}
\newcommand{\nT}{{\lfloor nT\rfloor}}
\newcommand{\proofend}{\hfill\mbox{$\Box$}}

\numberwithin{equation}{section}

\theoremstyle{change} \theorembodyfont{\em}
\newtheorem{Lem}{Lemma.}[section]
\newtheorem{Thm}[Lem]{Theorem.}
\newtheorem{Pro}[Lem]{Proposition.}
\newtheorem{Cor}[Lem]{Corollary.}
\newtheorem{Def}[Lem]{Definition.}

\theorembodyfont{\rm}
\newtheorem{Rem}[Lem]{Remark.}
\newtheorem{Ex}[Lem]{Example.}

\begin{document}

\begin{center}
 {\bfseries\Large Asymptotic behavior of critical, irreducible \\
                   multi-type continuous state and continuous time \\[3mm]
                   branching processes with immigration}
 \\[6mm]

 {\sc\large
  M\'aty\'as $\text{Barczy}^{*,\diamond}$,
  \ Gyula $\text{Pap}^{**}$}

\end{center}

\vskip0.2cm

\noindent
 * Faculty of Informatics, University of Debrecen,
   Pf.~12, H--4010 Debrecen, Hungary.

\noindent
 ** Bolyai Institute, University of Szeged,
     Aradi v\'ertan\'uk tere 1, H--6720 Szeged, Hungary.

\noindent e--mails: barczy.matyas@inf.unideb.hu (M. Barczy),
                    papgy@math.u-szeged.hu (G. Pap).

\noindent $\diamond$ Corresponding author.

\vskip0.2cm


\renewcommand{\thefootnote}{}
\vspace*{0.2cm}
\footnote{The research was realized in the frames of
 T\'AMOP 4.2.4.\ A/2-11-1-2012-0001 ,,National Excellence Program --
 Elaborating and operating an inland student and researcher personal support
 system''.
The project was subsidized by the European Union and co-financed by the
 European Social Fund.}

\vspace*{-10mm}

\begin{abstract}
Under natural assumptions, a Feller type diffusion approximation is derived for
 critical, irreducible multi-type continuous state and continuous time
 branching processes with immigration.
Namely, it is proved that a sequence of appropriately scaled random step
 functions formed from a critical, irreducible multi-type continuous state
 and continuous time branching process with immigration converges weakly
 towards a squared Bessel process supported by a ray determined by the Perron
 vector of a matrix related to the branching mechanism of the branching
 process in question.
\end{abstract}

\section{Introduction}
\label{section_intro}

The study of the limit behaviour of critical multi-type branching processes has a
 long tradition and history.
Most of the literature is devoted to so-called conditioned limit theorems for branching processes without
 immigration.

For a certain class of critical positively regular discrete time branching processes
 without immigration \ $(Z_n)_{n\geq1}$ \ with a finite or infinite number of types,
 Mullikin \cite[Theorem 9]{Mul} characterized the conditional limiting distribution
 of \ $n^{-1} Z_n$ \ given that \ $Z_n \ne 0$ \ as \ $n \to \infty$.

For critical discrete time branching processes with immigration \ $(X_n)_{n\geq 1}$, \
 under second order moment assumptions, Wei and Winnicki \cite[Theorem 2.1]{WeiWin} proved
 weak convergence of a sequence of step processes \ $(n^{-1}X_{\lfloor nt\rfloor})_{t\geq0}$, $n\geq 1$,
 \ as \ $n\to\infty$, \ characterizing the limit process as a squared Bessel process.

Ren et al.~\cite{RenYanZha} investigated conditional limit theorems for critical
 continuous-state and continuous time branching processes \ $(X_t)_{t\geq0}$ \ with branching mechanism \ $\lambda^{1+\alpha} L(1/\lambda)$, \
 but without immigration, where \ $\alpha\in[0,1]$ \ and \ $L$ \ is slowly varying at \ $\infty$.
\ They proved that if \ $\alpha\in(0, 1]$, \ then there are norming constants \ $Q_t\to0$ \ (as \ $t\uparrow \infty$)
 \ such that for every \ $x > 0$, \ $P_x(Q_tX_t \in\cdot \mid X_t > 0)$ \ converges weakly to a non-degenerate limit
 \ as \ $t\to\infty$.
\ As a continuation of these investigations, Ren et al.~\cite{RenSonZha} studied
 conditional limit theorems for some critical superprocesses conditioned on non-extinction.

Recently, Iyer et al.~\cite[Section 7]{IyeLegPeg} investigated limit theorems for critical
 continuous time and continuous state branching processes (without immigration)
 that become extinct almost surely.
First, they obtained a conditional limit theorem for fixed initial population size,
 and next, they studied non-conditioned scaling limits with initial
 population size scaled to obtain non-degenerate L\'evy process limits.

For a single-type (one-dimensional) critical continuous state and continuous
 time branching process with immigration (CBI process) \ $X$, \ under second order
 moment assumptions, Huang et al.~\cite[Theorem 2.3]{HuaMaZhu}
 characterized the limiting distribution of
 \ $(n^{-1} X_{\lfloor nt\rfloor})_{t\geq 0}$ \ as \ $n\to\infty$.
\ Our Theorem \ref{main} is a generalization of this result to the multi-type case (under fourth order moment assumptions),
 see Remark \ref{REMARK2}, and it may serve as a key tool for studying asymptotic behaviour of conditional least squares estimators
 of some parameters of processes in question.

The paper is organized as follows.
In Section \ref{section_CBI}, for completeness and better readability, we
 recall from Barczy et al.~\cite{BarLiPap2} some notions and statements for
 multi-type CBI processes such as the form of their infinitesimal generator,
 and a formula for their first moment.
In Section \ref{section_classification}, we introduce irreducible CBI processes and
 we give a classification, namely we define subcritical, critical and supercritical
 irreducible CBI processes, see Definitions \ref{Def_irreducible} and \ref{Def_indecomposable_crit},
 respectively.

In Section \ref{Convergence_result}, it is proved that the sequence
 \ $(n^{-1} \bX_{\lfloor nt\rfloor})_{t\geq 0}$, $n\geq 1$, \
 of scaled random step functions formed from a critical, irreducible multi-type CBI process \ $\bX$ \
 converges weakly towards a squared Bessel process supported by a ray determined by the Perron
 vector of a matrix related to the branching mechanism of \ $\bX$, \ see Theorem \ref{main}.
The limit process is characterized as a pathwise unique strong solution of a stochastic differential equation (SDE).
In Remark \ref{REMARK2}, we specialize Theorem \ref{main} to dimension 1 noting that the limit process
 is a single-type CBI diffusion process.
In Remark \ref{REMARK3} we point out that in case of \ $d\geq 2$, \ the limit process is not a
 \ $d$-type CBI process any longer, contrary to dimension 1.
We also formulate a consequence of Theorem \ref{main} deriving a limit distribution for the
 relative frequencies of distinct types of individuals, see Corollary \ref{Cor_rel_freq}.
For different models, one can find similar results in Jagers \cite[Corollary 1]{Jag}
 and in Yakovlev and Yanev \cite[Theorem 2]{YakYan}.
Section \ref{Examples} is devoted to give examples for multi-type CBI processes so that
 the drift and diffusion coefficients of the SDE characterizing the limit process in Theorem \ref{main}
 are calculated explicitly.
In Section \ref{Proof_main} we prove Theorem \ref{main}.
First, we prove weak convergence of a sequence of martingale differences \ $(\bM_n)_{n\geq 1}$ \ constructed from
 an irreducible and critical multi-type CBI process \ $\bX$.
Namely, \ $\bM_n$ \ is the difference of \ $\bX_n$ \ and the conditional expectation of \ $\bX_n$ \
 with respect to \ $\bX_{n-1}$.
\ The proof is based on a result due to Isp\'any and Pap \cite[Corollary 2.2]{IspPap2}
 (see also Theorem \ref{Conv2DiffThm}), which is about convergence of
 random step processes towards a diffusion process.
Using weak convergence of \ $(\bM_n)_{n\geq 1}$, \ an application of a version of the
 continuous mapping theorem (see Lemma \ref{Conv2Funct}) yields weak convergence of
 \ $(n^{-1} \bX_{\lfloor nt\rfloor})_{t\geq 0}$ \ as \ $n\to\infty$.
Comparing our proof of technique with that of
 Huang et al.~\cite[Theorem 2.3]{HuaMaZhu} (where the dimension is \ $1$), \
 they are completely different.
Huang et al.~\cite{HuaMaZhu} start with a SDE with jumps for the process \ $X$,
 \ and, applying Doob's inequality, tightness of the scaled processes
 \ $(n^{-1} X_{nt})_{t\geq 0}$, \ $n\geq 1$, \ is shown.
Then, by Skorokhod's theorem and a semimartingale representation theorem,
 they prove weak convergence of \ $(n^{-1} X_{nt})_{t\geq 0}$ \ as \ $n\to\infty$
 \ characterizing the limit distribution as well.
Finally, an application of the continuous mapping theorem yields weak convergence of
 \ $(n^{-1} X_{\lfloor nt\rfloor})_{t\geq 0}$ \ as \ $n\to\infty$ \ with the same limit distribution.
Comparing our technique of proof with that of Wei and Winnicki \cite[Theorem 2.1]{WeiWin},
 one can realize that they are completely different as well.
They calculated the infinitesimal generator of step processes \ $(n^{-1}X_{\lfloor nt\rfloor})_{t\geq0}$, \ $n\geq 1$,
 \ and examined its limit behaviour as \ $n\to\infty$ \ resulting the infinitesimal
 generator of the limit squared Bessel process.
In a companion paper Barczy and Pap \cite[Corollary 3.5 and Remark 3.6]{BarPap} we investigate
 convergence properties of the sequence of infinitesimal generators of
 \ $(n^{-1}\bX_{\lfloor nt\rfloor})_{t\geq0}$, \ $n \geq 1$.
\ It is an open question whether our main Theorem \ref{main} might be proved by the help of
 infinitesimal generators.
Further, we note that, to the best knowledge of the authors, it is not known, whether
 the sequence of scaled processes \ $(n^{-1}\bX_{nt})_{t\geq0}$, \ $n \geq 1$, \ is weakly
 convergent for an irreducible and critical $d$-type CBI process \ $\bX$ \ with \ $d \geq 2$.

In Appendix \ref{section_moments} we present some Frobenius--Perron type results for irreducible matrices
 having non-negative off-diagonal entries.
Appendix \ref{section_moments_CBI} is devoted to study asymptotic behaviour of moments of
 irreducible and critical multi-type CBI processes.
First, we describe the asymptotic behaviour of the first moment of irreducible multi-type CBI processes, see Proposition \ref{mean_asymptotics}.
The classification given in Definition \ref{Def_indecomposable_crit} is based
 on this description.
In case of an irreducible and critical multi-type CBI process \ $\bX$,
 \ we derive some moment estimations for the process and the corresponding sequence of martingale differences
 \ $(\bM_n)_{n\geq 1}$, \ see Lemmas \ref{moment_estimations_X_critical} and \ref{moment_estimations_1_2}, respectively.
In Appendix \ref{section_conv_step_drocesses}, we recall a result about convergence of random step processes towards a
 diffusion process due to Isp\'any and Pap \cite[Corollary 2.2]{IspPap2}.
In Appendix \ref{continuous_mapping_theorem} we present a version of the continuous mapping theorem.

\section{Multi-type CBI processes}
\label{section_CBI}

Let \ $\ZZ_+$, \ $\NN$, \ $\RR$, \ $\RR_+$  \ and \ $\RR_{++}$ \ denote the set
 of non-negative integers, positive integers, real numbers, non-negative real
 numbers and positive real numbers, respectively.
For \ $x , y \in \RR$, \ we will use the notations
 \ $x \land y := \min \{x, y\} $ \ and \ $x^+:= \max \{0, x\} $.
\ By \ $\|\bx\|$ \ and \ $\|\bA\|$, \ we denote the Euclidean norm of a vector
 \ $\bx \in \RR^d$ \ and the induced matrix norm of a matrix
 \ $\bA \in \RR^{d\times d}$, \ respectively.
The natural basis in \ $\RR^d$ \ and the Borel $\sigma$-algebras on \ $\RR^d$
 \ and on \ $\RR_+^d$ \ will be denoted by \ $\be_1$, \ldots, $\be_d$, \ and
 by \ $\cB(\RR^d)$ \ and \ $\cB(\RR_+^d)$, \ respectively.
The $d$-dimensional unit matrix is denoted by \ $\bI_d$.
\ For \ $\bx = (x_i)_{i\in\{1,\ldots,d\}} \in \RR^d$ \ and
 \ $\by = (y_i)_{i\in\{1,\ldots,d\}} \in \RR^d$, \ we will use the notation
 \ $\bx \leq \by$ \ indicating that \ $x_i \leq y_i$ \ for all
 \ $i \in \{1, \ldots, d\}$.
\ By \ $C^2_\cc(\RR_+^d,\RR)$ \ we denote the set of twice continuously
 differentiable real-valued functions on \ $\RR_+^d$ \ with compact support.
Throughout this paper, we make the conventions \ $\int_a^b := \int_{(a,b]}$
 \ and \ $\int_a^\infty := \int_{(a,\infty)}$ \ for any \ $a, b \in \RR$ \ with
 \ $a < b$.

\begin{Def}\label{Def_essentially_non-negative}
A matrix \ $\bA = (a_{i,j})_{i,j\in\{1,\ldots,d\}} \in \RR^{d\times d}$ \ is called
 essentially non-negative if \ $a_{i,j} \in \RR_+$ \ whenever
 \ $i, j \in \{1,\ldots,d\}$ \ with \ $i \ne j$, \ i.e., if \ $\bA$ \ has
 non-negative off-diagonal entries.
The set of essentially non-negative \ $d \times d$ \ matrices will be denoted
 by \ $\RR^{d\times d}_{(+)}$.
\end{Def}

\begin{Def}\label{Def_admissible}
A tuple \ $(d, \bc, \Bbeta, \bB, \nu, \bmu)$ \ is called a set of admissible
 parameters if
 \renewcommand{\labelenumi}{{\rm(\roman{enumi})}}
 \begin{enumerate}
  \item
   $d \in \NN$,
  \item
   $\bc = (c_i)_{i\in\{1,\ldots,d\}} \in \RR_+^d$,
  \item
   $\Bbeta = (\beta_i)_{i\in\{1,\ldots,d\}} \in \RR_+^d$,
  \item
   $\bB = (b_{i,j})_{i,j\in\{1,\ldots,d\}} \in \RR^{d \times d}_{(+)}$,
  \item
   $\nu$ \ is a Borel measure on \ $U_d := \RR_+^d \setminus \{\bzero\}$
    \ satisfying \ $\int_{U_d} (1 \land \|\bz\|) \, \nu(\dd \bz) < \infty$,
  \item
   $\bmu = (\mu_1, \ldots, \mu_d)$, \ where, for each
    \ $i \in \{1, \ldots, d\}$, \ $\mu_i$ \ is a Borel measure on \ $U_d$
    \ satisfying
    \begin{align}\label{help2_intcond_mu}
      \int_{U_d}
       \left[ \|\bz\| \wedge \|\bz\|^2
              + \sum_{j \in \{1, \ldots, d\} \setminus \{i\}} z_j \right]
       \mu_i(\dd \bz)
      < \infty .
    \end{align}
  \end{enumerate}
\end{Def}

\begin{Rem}
Our Definition \ref{Def_admissible} of the set of admissible
 parameters is a special case of Definition 2.6 in
 Duffie et al.~\cite{DufFilSch}, which is suitable for all affine processes,
 see Barczy et al.~\cite[Remark 2.3]{BarLiPap2}.
Further, for all \ $i \in \{1, \ldots, d\}$, \ condition \eqref{help2_intcond_mu}
 is equivalent to
 \begin{align*}
  \int_{U_d}
     \left[ (1 \land z_i)^2
             + \sum_{j \in \{1, \ldots, d\} \setminus \{i\}} (1 \land z_j) \right]
      \mu_i(\dd \bz)
    < \infty
   \quad \text{and} \quad
   \int_{U_d} \|\bz\| \bbone_{\{\|\bz\|\geq 1\}}\,\mu_i(\dd \bz) < \infty ,
 \end{align*}
 see Barczy et al.~\cite[Remark 2.3]{BarLiPap2}.
\proofend
\end{Rem}

\begin{Thm}\label{CBI_exists}
Let \ $(d, \bc, \Bbeta, \bB, \nu, \bmu)$ \ be a set of admissible parameters.
Then there exists a unique conservative transition semigroup \ $(P_t)_{t\in\RR_+}$ \ acting on
 the Banach space (endowed with the supremum norm) of real-valued bounded
 Borel-measurable functions on the state space \ $\RR_+^d$ \ such that its
 infinitesimal generator is
 \begin{equation}\label{CBI_inf_gen}
  \begin{aligned}
   (\cA f)(\bx)
   &= \sum_{i=1}^d c_i x_i f_{i,i}''(\bx)
      + \langle \Bbeta + \bB \bx, \Bf'(\bx) \rangle
      + \int_{U_d} \bigl( f(\bx + \bz) - f(\bx) \bigr) \, \nu(\dd \bz) \\
   &\phantom{\quad}
      + \sum_{i=1}^d
         x_i
         \int_{U_d}
          \bigl( f(\bx + \bz) - f(\bx) - f'_i(\bx) (1 \land z_i) \bigr)
          \, \mu_i(\dd \bz)
  \end{aligned}
 \end{equation}
 for \ $f \in C^2_\cc(\RR_+^d,\RR)$ \ and \ $\bx \in \RR_+^d$, \ where \ $f_i'$
\ and \ $f_{i,i}''$, \ $i \in \{1, \ldots, d\}$, \ denote the first and second
 order partial derivatives of \ $f$ \ with respect to its \ $i$-th variable,
 respectively, and \ $\Bf'(\bx) := (f_1'(\bx), \ldots, f_d'(\bx))^\top$.
\end{Thm}

\begin{Rem}
This theorem is a special case of Theorem 2.7 of Duffie et
 al.~\cite{DufFilSch} with \ $m = d$, \ $n = 0$ \ and zero killing rate.
\proofend
\end{Rem}

\begin{Def}\label{Def_CBI}
A conservative Markov process with state space \ $\RR_+^d$ \ and with transition semigroup
 \ $(P_t)_{t\in\RR_+}$ \ given in Theorem \ref{CBI_exists} is called a multi-type
 CBI process with parameters \ $(d, \bc, \Bbeta, \bB, \nu, \bmu)$.
\end{Def}

Let \ $(\bX_t)_{t\in\RR_+}$ \ be a multi-type CBI process with parameters
 \ $(d, \bc, \Bbeta, \bB, \nu, \bmu)$ \ such that \ $\EE(\|\bX_0\|) < \infty$
 \ and the moment condition
 \begin{equation}\label{moment_condition_1}
  \int_{U_d} \|\bz\| \bbone_{\{\|\bz\|\geq1\}} \, \nu(\dd \bz) < \infty
 \end{equation}
 holds.
Then, by Lemma 3.4 in Barczy et al. \cite{BarLiPap2},
 \begin{equation}\label{EXbX}
  \EE(\bX_t) = \ee^{t\tbB} \EE(\bX_0)
               + \left( \int_0^t \ee^{u\tbB} \, \dd u \right) \tBbeta ,
  \qquad t \in \RR_+ ,
 \end{equation}
 where
 \begin{gather}
  \tbB := (\tb_{i,j})_{i,j\in\{1,\ldots,d\}} , \qquad
  \tb_{i,j} := b_{i,j}
              + \int_{U_d} (z_i - \delta_{i,j})^+ \, \mu_j(\dd \bz) ,
  \label{tbB} \\
  \tBbeta := \Bbeta + \int_{U_d} \bz \, \nu(\dd \bz) ,
  \label{tBbeta}
 \end{gather}
 with \ $\delta_{i,j}:=1$ \ if \ $i=j$, \ and \ $\delta_{i,j}:=0$ \ if
 \ $i \ne j$.
\ Note that \ $\tbB \in \RR^{d \times d}_{(+)}$ \ and \ $\tBbeta \in \RR_+^d$,
 \ since
 \begin{equation}\label{help}
  \int_{U_d} \|\bz\| \, \nu(\dd\bz) < \infty , \qquad
  \int_{U_d} (z_i - \delta_{i,j})^+ \, \mu_j(\dd \bz) < \infty , \quad
  i, j \in \{1, \ldots, d\} ,
 \end{equation}
 see Barczy et al. \cite[Section 2]{BarLiPap2}.

\section{Classification of multi-type CBI processes and
             moment estimations}
\label{section_classification}

For a matrix \ $\bA \in \RR^{d \times d}$, \ $\sigma(\bA)$ \ will denote the
 spectrum of \ $\bA$, \ i.e., the set of the eigenvalues of \ $\bA$.
\ Then \ $r(\bA) := \max_{\lambda \in \sigma(\bA)} |\lambda|$ \ is the spectral
 radius of \ $\bA$.
\ Moreover, we will use the notation
 \[
   s(\bA) := \max_{\lambda \in \sigma(\bA)} \Re(\lambda) .
 \]
By the spectral mapping theorem (see, e.g., Dunford and Schwartz
 \cite[Theorem VII.3.11]{DunSch}), \ $\sigma(\ee^{t\bA}) = \ee^{t\sigma(\bA)}$
 \ for all \ $t \in \RR_+$.
\ Consequently,
 \begin{equation}\label{rs}
  r(\ee^{t\bA}) = \max_{\lambda \in \sigma(\bA)} |\ee^{t\lambda}|
  = \max_{\lambda \in \sigma(\bA)} \ee^{t\Re(\lambda)}
  = \ee^{s(\bA)t} , \qquad t \in \RR_+ ,
 \end{equation}
  and hence \ $s(\bA) = \log r(\ee^\bA)$.
\ A matrix \ $\bA \in \RR^{d\times d}$ \ is called reducible if there exist a
 permutation matrix \ $\bP \in \RR^{ d \times d}$ \ and an integer \ $r$ \ with
 \ $1 \leq r \leq d-1$ \ such that
 \[
  \bP^\top \bA \bP
   = \begin{bmatrix} \bA_1 & \bA_2 \\ \bzero & \bA_3 \end{bmatrix},
 \]
 where \ $\bA_1 \in \RR^{r\times r}$, \ $\bA_3 \in \RR^{ (d-r) \times (d-r) }$,
 \ $\bA_2 \in \RR^{ r \times (d-r) }$, \ and \ $\bzero \in \RR^{ (d-r) \times r}$
 \ is a null matrix.
A matrix \ $\bA \in \RR^{d\times d}$ \ is called irreducible if it is not
 reducible, see, e.g.,
 Horn and Johnson \cite[Definitions 6.2.21 and 6.2.22]{HorJoh}.
We do emphasize that no 1-by-1 matrix is reducible.

If \ $(\bX_t)_{t\in\RR_+}$ \ is a CBI process with parameters
 \ $(d, \bc, \Bbeta, \bB, \nu, \bmu)$ \ such that the moment condition
 \eqref{moment_condition_1} holds, then \ $\ee^{t\tbB} \in \RR_+^{d \times d}$
 \ for all \ $t \in \RR_+$, \ since \ $\tbB \in \RR^{d \times d}_{(+)}$, \ see the
 explanation before Lemma \ref{FP1}.
Moreover, by Lemma \ref{FP1} and Remark \ref{FP1_Rem},
 \ $\ee^{t_0\tbB} \in \RR_{++}^{d \times d}$ \ for some (and hence for all)
 \ $t_0 \in \RR_{++}$ \ if and only if \ $\tbB$ \ is irreducible or,
 if and only if \ $\ee^{\tbB}$ \ is irreducible.

\begin{Def}\label{Def_irreducible}
Let \ $(\bX_t)_{t\in\RR_+}$ \ be a multi-type CBI process with parameters
 \ $(d, \bc, \Bbeta, \bB, \nu, \bmu)$ \ such that the moment condition
 \eqref{moment_condition_1} holds.
Then \ $(\bX_t)_{t\in\RR_+}$ \ is called irreducible if \ $\tbB$ \ is
 irreducible.
\end{Def}

Next we introduce a classification of irreducible multi-type CBI processes.
Formula \eqref{EXbX} shows that the semigroup
 \ $\bigl( \ee^{t\tbB} \bigr)_{t\in\RR_+}$ \ of matrices plays a crucial role in
 the asymptotic behavior of the expectations \ $\EE(\bX_t)$ \ as
 \ $t \to \infty$ \ described in Proposition \ref{mean_asymptotics}.
This gives a motivation for a classification of irreducible multi-type CBI
 processes.

\begin{Def}\label{Def_indecomposable_crit}
Let \ $(\bX_t)_{t\in\RR_+}$ \ be a multi-type CBI process with parameters
 \ $(d, \bc, \Bbeta, \bB, \nu, \bmu)$ \ such that \ $\EE(\|\bX_0\|) < \infty$
 \ and the moment condition \eqref{moment_condition_1} holds.
Suppose that \ $(\bX_t)_{t\in\RR_+}$ \ is irreducible.
Then \ $(\bX_t)_{t\in\RR_+}$ \ is called
 \[
   \begin{cases}
    subcritical & \text{if \ $s(\tbB) < 0$,} \\
    critical & \text{if \ $s(\tbB) = 0$,} \\
    supercritical & \text{if \ $s(\tbB) > 0$.}
   \end{cases}
 \]
\end{Def}

The classification for subcritical, critical and supercritical cases in
 Definition \ref{Def_indecomposable_crit} is in accordance with the
 corresponding classification for single-type continuous state and continuous time branching processes,
 see, e.g., Li \cite[page 58]{Li}.

\section{Convergence result}
\label{Convergence_result}

A function \ $f : \RR_+ \to \RR^d$ \ is called \emph{c\`adl\`ag} if it is right
 continuous with left limits.
\ Let \ $\DD(\RR_+, \RR^d)$ \ and \ $\CC(\RR_+, \RR^d)$ \ denote the space of
 all \ $\RR^d$-valued c\`adl\`ag and continuous functions on \ $\RR_+$,
 \ respectively.
Let \ $\cD_\infty(\RR_+, \RR^d)$ \ denote the Borel $\sigma$-field in
 \ $\DD(\RR_+, \RR^d)$ \ for the metric characterized by Jacod and Shiryaev
 \cite[VI.1.15]{JacShi} (with this metric, \ $\DD(\RR_+, \RR^d)$ \ is a complete
 and separable metric space).
For \ $\RR^d$-valued stochastic processes \ $(\bcY_t)_{t\in\RR_+}$ \ and
 \ $(\bcY^n_t)_{t\in\RR_+}$, \ $n \in \NN$, \ with c\`adl\`ag paths we write
 \ $\bcY^n \distr \bcY$ \ as \ $n\to\infty$ \ if the distribution of
 \ $\bcY^n$ \ on the space \ $(\DD(\RR_+, \RR^d), \cD_\infty(\RR_+, \RR^d))$
 \ converges weakly to the distribution of \ $\bcY$ \ on the space
 \ $(\DD(\RR_+, \RR^d), \cD_\infty(\RR_+, \RR^d))$ \ as \ $n \to \infty$.

\begin{Thm}\label{main}
Let \ $(\bX_t)_{t\in\RR_+}$ \ be a multi-type CBI process with
 parameters \ $(d, \bc, \Bbeta, \bB, \nu, \bmu)$ \ such that
 \ $\EE(\|\bX_0\|^4) < \infty$ \ and
 \begin{equation}\label{moment_condition_4}
  \int_{U_d} \|\bz\|^4 \bbone_{\{\|\bz\|\geq1\}} \, \nu(\dd \bz) < \infty , \qquad
  \int_{U_d} \|\bz\|^4 \bbone_{\{\|\bz\|\geq1\}} \, \mu_i(\dd \bz) < \infty , \quad
  i \in \{1, \ldots, d\} .
 \end{equation}
Suppose that \ $(\bX_t)_{t\in\RR_+}$ \ is irreducible and critical.
Then
 \begin{gather}\label{Conv_X}
  (\bcX_t^{(n)})_{t\in\RR_+} := (n^{-1} \bX_{\nt})_{t\in\RR_+}
  \distr (\bcX_t)_{t\in\RR_+} := (\cX_t \bu)_{t\in\RR_+} \qquad
  \text{as \ $n \to \infty$}
 \end{gather}
 in \ $\DD(\RR_+, \RR^d)$, \ where \ $\bu := \bu_{\tbB} \in \RR_{++}^d$ \ is the
 right Perron vector of \ $\ee^{\tbB}$ \ corresponding to the eigenvalue \ $1$
 \ with \ $\sum_{i=1}^d \be_i^\top \bu = 1$ \ (see \textup{(ii)} of Lemma \ref{FP}),
 \ $(\cX_t)_{t \in \RR_+}$ \ is the unique strong solution of the SDE
 \begin{equation}\label{SDE_X}
  \dd \cX_t
  = \langle \bv, \tBbeta\rangle \, \dd t
    + \sqrt{\langle \obC \bv, \bv \rangle \cX_t^+} \, \dd \cW_t , \qquad t \in \RR_+ ,
  \qquad \cX_0 = 0 ,
 \end{equation}
 where \ $\bv := \bv_{\tbB} \in \RR_{++}^d$ \ is the left Perron vector of
 \ $\ee^{\tbB}$ \ corresponding to the eigenvalue \ $1$ \ with \ $\bv^\top \bu = 1$
 \ (see \textup{(iii)} of Lemma \ref{FP}), \ $(\cW_t)_{t \in \RR_+}$ \ is a standard
 Brownian motion, \ $\tBbeta$ \ is given in \eqref{tBbeta}, and
 \begin{gather}
  \obC := \sum_{k=1}^d \langle \be_k, \bu \rangle \bC_k  \in \RR_+^{d \times d}  \label{obC}
  \end{gather}
  with
  \begin{align}
  \bC_k := 2 c_k \be_k \be_k^\top + \int_{U_d} \bz \bz^\top \mu_k(\dd \bz)
        \in \RR_+^{d \times d} , \qquad
  k \in \{1, \ldots, d\} . \label{tbC}
 \end{align}
\end{Thm}

\begin{Rem}
We suspect that the moment conditions might be relaxed to
 \ $\EE(\|\bX_0\|^2) < \infty$ \ and
 \begin{equation}\label{moment_condition_2}
  \int_{U_d} \|\bz\|^2 \bbone_{\{\|\bz\|\geq1\}} \, \nu(\dd \bz) < \infty , \qquad
  \int_{U_d} \|\bz\|^2 \bbone_{\{\|\bz\|\geq1\}} \, \mu_i(\dd \bz) < \infty , \quad
  i \in \{1, \ldots, d\} .
 \end{equation}
In fact, the higher order moment assumptions are used only for checking the
 conditional Lindeberg condition, namely, condition (ii) of Theorem
 \ref{Conv2DiffThm}, in order to prove convergence \eqref{Conv_M} of an
 appropriately defined sequence of martingale differences \eqref{cMk}.
One might check the conditional Lindeberg condition under
 \ $\EE(\|\bX_0\|^2) < \infty$ \ and the above weaker moment assumptions
 \eqref{moment_condition_2} by the method of Isp\'any and Pap \cite{IspPap1},
 see also this method in Barczy et al.~\cite[proof of convergence (5.2)]{BarIspPap0}.
One might generalize Theorem \ref{main} to an appropriate sequence of initial distributions
 instead of a fixed one, see, e.g., the method of Isp\'any and Pap \cite{IspPap3}.
\proofend
\end{Rem}

\begin{Rem}\label{REMARK0}
Among the moment conditions we have the relationships
 \eqref{moment_condition_4} $\Rightarrow$ \eqref{moment_condition_2}
 $\Rightarrow$ \eqref{moment_condition_1}.
The moment conditions \eqref{moment_condition_2} together with the fact that
 \ $\nu$ \ and \ $\bmu$ \ satisfy Definition \ref{Def_admissible} imply
 \begin{equation}\label{moment_condition_2_m}
  \int_{U_d} \|\bz\|^2 \, \nu(\dd \bz) < \infty , \qquad
  \int_{U_d} \|\bz\|^2 \, \mu_k(\dd \bz) < \infty , \quad
  k \in \{1, \ldots, d\} .
 \end{equation}
Indeed,
 \begin{align*}
  \int_{U_d} \|\bz\|^2 \, \nu(\dd \bz)
  &= \int_{U_d} \|\bz\|^2 \bbone_{\{\|\bz\|< 1\}} \, \nu(\dd \bz)
     + \int_{U_d} \|\bz\|^2 \bbone_{\{\|\bz\|\geq1\}} \, \nu(\dd \bz) \\
  &\leq \int_{U_d} (1 \land \|\bz\|) \, \nu(\dd \bz)
        + \int_{U_d} \|\bz\|^2 \bbone_{\{\|\bz\|\geq1\}} \, \nu(\dd \bz)
   < \infty ,
 \end{align*}
 and for all \ $k \in \{1, \ldots, d\}$,
 \begin{align*}
  \int_{U_d} \|\bz\|^2 \, \mu_k(\dd \bz)
  &= \int_{U_d} \|\bz\|^2 \bbone_{\{\|\bz\|<1\}} \, \mu_k(\dd \bz)
      + \int_{U_d} \|\bz\|^2 \bbone_{\{\|\bz\|\geq1\}} \, \mu_k(\dd \bz) \\
  &\leq \int_{U_d} \|\bz\| \land \|\bz\|^2 \, \mu_k(\dd \bz)
        + \int_{U_d} \|\bz\|^2 \bbone_{\{\|\bz\|\geq1\}} \, \mu_k(\dd \bz)
   < \infty .
 \end{align*}
Clearly, \eqref{moment_condition_2_m} implies also
 \ $\bC_k \in \RR_+^{d \times d}$, \ $k \in \{1, \ldots, d\}$, \ since
 \[
   \int_{U_d} \|\bz \bz^\top\| \, \mu_k(\dd \bz)
   \leq \int_{U_d} \|\bz\|^2 \, \mu_k(\dd \bz)
   < \infty .
 \]
Obviously, \ $\bC_k$, \ $k \in \{1, \ldots, d\}$, \ and
 \ $\obC$ \ are symmetric positive semidefinite matrices,
 and \ $\obC = \bzero$ \ if and only if \ $\bC_k = \bzero$ \ for all
 \ $k \in \{1, \ldots, d\}$.
\ Indeed, \ $\obC = \bzero$ \ implies
 \begin{align*}
  0= \langle \obC\be_i,\be_i\rangle
   = \sum_{k=1}^d \be_k^\top \bu \langle \bC_k\be_i,\be_i\rangle,
   \qquad i\in\{1,\ldots,d\},
 \end{align*}
 and, since \ $\be_k^\top \bu\in\RR_{++}$ \ and \ $\langle \bC_k\be_i,\be_i\rangle\in\RR_+$ \
 (due to positive semidefiniteness of \ $\bC_k$),
 \ we get \ $\langle \bC_k\be_i,\be_i\rangle=0$, \ $k,i\in\{1,\ldots,d\}$.
\ For each \ $k,i\in\{1,\ldots,d\}$, \ we have
 \ $0 = \langle \bC_k\be_i,\be_i\rangle
      = \langle \sqrt{\bC_k} \be_i, \sqrt{\bC_k} \be_i \rangle = \|\sqrt{\bC_k} \be_i\|^2$,
 \ where \ $\sqrt{\bC_k}$ \ denotes the unique symmetric and positive semidefinite square root of
 \ $\bC_k$.
\ Consequently, \ $\sqrt{\bC_k} \be_i = \bzero$, \ thus
 \ $\bC_k \be_i = \sqrt{\bC_k} \sqrt{\bC_k} \be_i = \bzero$, \ implying
 \ $\langle \bC_k\be_i,\be_j\rangle = 0$ \ for each \ $k,i,j\in\{1,\ldots,d\}$.
\ Hence \ $\obC = \bzero$ \ if and only if \ $c_k = 0$ \ and
 \ $\mu_k = 0$ \ for all \ $k \in \{1, \ldots, d\}$.
\proofend
\end{Rem}

\begin{Rem}\label{REMARK1}
The SDE \eqref{SDE_X} has a pathwise unique strong solution
 \ $(\cX_t^{(x)})_{t\in\RR_+}$ \ for all initial values
 \ $\cX_0^{(x)} = x \in \RR$, \ and if the initial value \ $x$ \ is
 nonnegative, then \ $\cX_t^{(x)}$ \ is nonnegative for all \ $t \in \RR_+$
 \ with probability one, since \ $\langle \bv, \tBbeta \rangle \in \RR_+$.
In fact, \ $(4 \langle \obC \bv, \bv \rangle^{-1} \cX_t^{(x)})_{t\in\RR_+}$ \ is a
 square of a \ $4 \langle \obC \bv, \bv \rangle^{-1} \langle \bv, \tBbeta \rangle$-dimensional
 Bessel process started at \ $x$ \ (see, e.g., Revuz and Yor \cite[Definitions XI.1.1]{RevYor}).
Moreover, if \ $\langle \bv, \tBbeta \rangle > 0$ \ then \ $\PP(\cX_t > 0) = 1$ \ for all \ $t > 0$,
 \ where \ $\cX_t =  \cX_t^{(0)}$, \ $t \in \RR_+$.
\ Indeed, we have \ $\cX_t = \langle \bv, \tBbeta\rangle t$, \ $t \in \RR_+$, \ if
\ $\langle \obC \bv, \bv \rangle = 0$, \ and \ $\cX_t$ \ has a gamma distribution with parameters
 \ $2 \langle \bv, \tBbeta\rangle / \langle \obC \bv, \bv \rangle$ \ and
 \ $2 / (\langle \obC \bv, \bv \rangle t)$ \ if \ $\langle \obC \bv, \bv \rangle \ne 0$ \ and \ $t > 0$.
\ For the proofs, see, e.g., Ikeda and Watanabe \cite[Chapter IV, Example 8.2]{IkeWat}.
Note that \ $\langle \bv, \tBbeta\rangle > 0$ \ if and only if \ $\tBbeta \ne \bzero$, \ which is
 equivalent to \ $\Bbeta \ne \bzero$ \ or \ $\nu \ne 0$.
\proofend
\end{Rem}

In the next remark we specialize Theorem \ref{main} for dimension 1.

\begin{Rem}\label{REMARK2}
If \ $(X_t)_{t\in\RR_+}$ \ is a single-type (1-dimensional and hence
 irreducible) critical CBI process with parameters
 \ $(1, c, \beta, b, \nu, \mu)$ \ satisfying assumptions of Theorem
 \ref{main}, then
 \ $\tb := b + \int_1^\infty (z-1) \, \mu(\dd z) = 0$ \ (due to criticality and
 using that \ $s(\tb) = \tb$ \ in dimension 1), and
 \ $(n^{-1} X_{\nt})_{t\in\RR_+} \distr (\cX_t)_{t\in\RR_+}$ \ as \ $n \to \infty$,
 \ where \ $(\cX_t)_{t\in\RR_+}$ \ is a pathwise unique strong solution of the
 SDE \ $\dd \cX_t = \tbeta \, \dd t + \sqrt{\oC \cX_t^+} \, \dd \cW_t$,
 \ $t \in \RR_+$, \ with initial value \ $\cX_0 = 0$ \ and with
 \ $\tbeta := \beta + \int_{U_1} z \, \nu(\dd z)$,
 \ $\oC := 2 c + \int_{U_1} z^2 \, \mu(\dd z)$.
\ Here, by \eqref{CBI_inf_gen}, the infinitesimal generator of
 \ $(X_t)_{t\in\RR_+}$ \ is
 \begin{align*}
  (\cA_X f)(x)
  &= c x f''(x) + (\beta + bx) f'(x)
     + \int_{U_1} (f(x + z) - f(x)) \, \nu(\dd z) \\
  &\quad
     + x \int_{U_1} (f(x + z) - f(x) - f'(x) (1 \land z)) \, \mu(\dd z) \\
  &= \frac{1}{2} \oC x f''(x) + (\tbeta + bx) f'(x)
     + \int_{U_1} (f(x + z) - f(x) - f'(x) z) \, \nu(\dd z) \\
  &\quad
     + x \int_{U_1}
          \left( f(x + z) - f(x) - f'(x) (1\land z)
                 - \frac{1}{2} f''(x) z^2 \right)
          \mu(\dd z) \\
  &= \frac{1}{2} \oC x f''(x) + \tbeta f'(x)
     + \int_{U_1} (f(x + z) - f(x) - f'(x) z) \, \nu(\dd z)  \\
  &\quad
     + x \int_{U_1}
          \left(f(x + z) - f(x) - f'(x) z - \frac{1}{2} f''(x) z^2 \right)
          \mu(\dd z)
 \end{align*}
 for \ $f \in \CC_\cc^2(\RR_+, \RR)$ \ and \ $x \in \RR_+$, \ where the last
 equality follows by
 \begin{align*}
  b + \int_{U_1} (z - (1\wedge z))\,\mu(\dd z)
  =  b + \int_1^\infty (z - 1)\,\mu(\dd z)
  = \tb
  = 0 .
 \end{align*}
Further, the limit process \ $(\cX_t)_{t\in\RR_+}$ \ is a single-type
 (1-dimensional) CBI diffusion process with parameters
 \ $\bigl(1, \frac{1}{2} \oC, \tbeta, 0, 0, 0\bigr)$, \ and its infinitesimal
 generator takes the form
 \begin{align*}
  (\cA_{\cX} f)(x)
    = \frac{1}{2} \oC x f''(x) + \tbeta f'(x)
 \end{align*}
 for \ $f \in \CC_\cc^2(\RR_+, \RR)$ \ and \ $x \in \RR_+$, \ see, e.g.,
 Karatzas and Shreve \cite[Section 5.1]{KarShr}.
Note that under the conditions
 \[
   \int_{U_1} z \, \nu(\dd z) < \infty , \qquad
   \int_{U_1} z^2 \, \mu(\dd z) < \infty ,
 \]
 a stronger statement, namely, a scaling limit theorem
 \ $(n^{-1} X_{nt})_{t\in\RR_+} \distr (\cX_t)_{t\in\RR_+}$ \ is also valid, see
 Barczy et al.~\cite[Corollary 2.1]{BarDorLiPap} (under the additional stronger
 assumption \ $\int_{U_1} z^2 \, \nu(\dd z) < \infty$, \
 Huang et al.~\cite[Theorem 2.3]{HuaMaZhu} proved this result with another method).
\proofend
\end{Rem}

In the next remark we point out that in case of \ $d \geq 2$, \ the limit
 process \ $(\bcX_t)_{t\in\RR_+}$ \ in \eqref{Conv_X} is not a \ $d$-type
 ($d$-dimensional) CBI process.

\begin{Rem}\label{REMARK3}
If \ $d \geq 2$ \ and \ $(\bX_t)_{t\in\RR_+}$ \ is an irreducible and critical
 \ $d$-type CBI process with parameters \ $(d, \bc, \Bbeta, \bB, \nu, \bmu)$
 \ satisfying assumptions of Theorem \ref{main}, then the limit process
 \ $(\bcX_t)_{t\in\RR_+}$ \ is not a \ $d$-type
 CBI process, i.e., in case \ $d \geq 2$, \ we have a different limit behaviour
 compared to dimension 1, see Remark \ref{REMARK2}.
Indeed, \ $(\bcX_t)_{t\in\RR_+}$ \ is a time homogeneous Markov process with
 state space \ $\RR_+ \bu$ \ and with infinitesimal generator
 \[
   (\cA_{\bcX} f)(\bx)
   = \frac{1}{2} \bv^\top \obC \bv x \tf^{''}(x) + \bv^\top \tBbeta \tf^{'}(x) ,
   \qquad \bx = x \bu , \qquad x \in \RR_+
 \]
 for \ $f : \RR_+ \bu \to \RR$, \ $f(x \bu) = \tf(x)$ \ with
 \ $\tf \in \CC_\cc^2(\RR_+, \RR)$, \ since
 \begin{align*}
  (\cA_{\bcX} f)(\bx)
  &= \lim_{h\downarrow0} h^{-1} [\EE(f(\bcX_h) \mid \bcX_0 = \bx) - f(\bx)]
   = \lim_{h\downarrow0} h^{-1}
      [\EE(f(\cX_h \bu) \mid \cX_0 \bu = x \bu) - f(x \bu)] \\
  &= \lim_{h\downarrow0} h^{-1} [\EE(\tf(\cX_h) \mid \cX_0 = x) - \tf(x)]
   = (\cA_\cX \tf)(x)
 \end{align*}
 has the above form by \eqref{SDE_X}, see, e.g., Karatzas and Shreve
 \cite[Section 5.1]{KarShr}.
Clearly, the infinitesimal generator \ $\cA_{\bcX}$ \ is not of the form
 \eqref{CBI_inf_gen}, \ since it is not defined for all
 \ $f \in \CC^2_\cc(\RR_+^d,\RR)$ \ due to \ $d \geq 2$.
\ Note that the process \ $(\cX_t)_{t\in\RR_+}$ \ is a single-type
 (one-dimensional) CBI process with parameters
 \ $\bigl(1, \frac{1}{2}\bv^\top \obC \bv, \bv^\top \tBbeta, 0, 0, 0\bigr)$.
\proofend
\end{Rem}

Next we formulate a corollary of Theorem \ref{main} deriving a limit distribution for the
 relative frequencies of distinct types of individuals.

\begin{Cor}\label{Cor_rel_freq}
Assume that the conditions of Theorem \ref{main} holds.
If, in addition, \ $\Bbeta \ne \bzero$ \ or \ $\nu \ne 0$, \ then for each \ $t > 0$ \ and
 \ $i, j \in \{1, \ldots, d\}$,
 \[
   \frac{\be_i^\top \bX_\nt}{\be_j^\top \bX_\nt} \stoch \frac{\be_i^\top \bu}{\be_j^\top \bu}
   \qquad \text{and} \qquad
   \frac{\be_i^\top \bX_\nt}{\sum_{k=1}^d \be_k^\top \bX_\nt} \stoch \be_i^\top \bu
   \qquad \text{as \ $n \to \infty$.}
 \]
\end{Cor}

\noindent
\textbf{Proof.}
Theorem \ref{main} implies
 \ $(\be_i^\top \bX_\nt, \be_j^\top \bX_\nt) \distr (\be_i^\top \bu \cX_t, \be_j^\top \bu \cX_t)$
 \ as \ $n \to \infty$ \ for each \ $t \in \RR_+$.
\ The function \ $g : \RR^2 \to \RR$, \ defined by
 \begin{align*}
  g(x, y)
  :=\begin{cases}
     \frac{x}{y} , & \text{if \ $x\in\RR$ \ and \ $y \ne 0$,}\\
     0 , & \text{if \ $x\in\RR$ \ and \ $y = 0$,}
    \end{cases}
 \end{align*}
 is continuous on the set \ $\RR\times(\RR \setminus \{0\})$, \ and the distribution of
 \ $(\be_i^\top \bu \cX_t, \be_j^\top \bu \cX_t)$
 \ is concentrated on this set, since, by Remark \ref{REMARK1}, \ $\PP(\cX_t > 0) = 1$, \ and
 \ $\be_j^\top \bu > 0$.
 \ Hence the continuous mapping theorem yields that
 \begin{align*}
  g(\be_i^\top \bX_\nt, \be_j^\top \bX_\nt)
  \distr
  g(\be_i^\top \bu \cX_t, \be_j^\top \bu \cX_t)
  = \frac{\be_i^\top \bu}{\be_j^\top \bu}  \qquad
  \text{as \ $n \to \infty$,}
 \end{align*}
 thus we obtain the first convergence.
Moreover,
 \[
   \frac{\be_i^\top \bX_\nt}{\sum_{k=1}^d \be_k^\top \bX_\nt}
   = \frac{1}{\sum_{k=1}^d \frac{\be_k^\top \bX_\nt}{\be_i^\top \bX_\nt}}
   \stoch
   \frac{1}{\sum_{k=1}^d \frac{\be_k^\top \bu}{\be_i^\top \bu}}
   = \frac{\be_i^\top \bu}{\sum_{k=1}^d \be_k^\top \bu}
   = \be_i^\top \bu , \qquad
   \text{as \ $n \to \infty$,}
 \]
 hence we obtain the second convergence.
\proofend

\section{Examples}
\label{Examples}

In this section we give some examples for multi-type CBI processes with
 parameters \ $(d, \bc, \Bbeta, \bB, \nu, \bmu)$ \ for which
 we give the vector \ $\tBbeta \in \RR_+^d$ \ and the matrix
 \ $\obC \in \RR_+^{d\times d}$, \ respectively, which appear in the drift and
 diffusion coefficients of the SDE \eqref{SDE_X}, respectively.

\begin{Ex} \
Let \ $(\bX_t)_{t\in\RR_+}$ \ be a 2-type (2-dimensional) CBI process with
 parameters \ $(2, \bc, \Bbeta, \bB, \nu, \bmu)$ \ such that
 \[
   \tbB = \gamma \begin{bmatrix} -1 & 1 \\ 1 & -1 \end{bmatrix}
 \]
 with some \ $\gamma \in \RR_{++}$.
\ Let us suppose that \ $\EE(\|\bX_0\|) < \infty$ \ and the moment condition
 \eqref{moment_condition_1} holds.
Then \ $\tbB$ \ is irreducible, and the eigenvalues of \ $\tbB$ \ are \ $0$
 \ and \ $-2\gamma$ \ thus \ $s(\tbB) = 0$, \ and hence \ $(\bX_t)_{t\in\RR_+}$
 \ is irreducible and critical.
Further, we have
 \begin{gather*}
  \bu = \bu_{\tbB} = \frac{1}{2} \begin{bmatrix} 1 \\ 1 \end{bmatrix} , \qquad
  \bv = \bv_{\tbB} = \begin{bmatrix} 1 \\ 1 \end{bmatrix} ,\qquad
  \obC = \frac{1}{2}(\bC_1+\bC_2).
 \end{gather*}
\end{Ex}

\begin{Ex} \
Let \ $(\bX_t)_{t\in\RR_+}$ \ be an irreducible and critical $d$-type
 ($d$-dimensional) CBI process with parameters
 \ $(d, \bzero, \bzero, \bzero, \nu, \bmu)$.
\ Let us suppose that the moment condition \eqref{moment_condition_1} holds.
Then
 \begin{gather*}
  \tbB = (\tb_{i,j})_{i,j\in\{1,\ldots,d\}}  \qquad \text{with} \qquad
  \tb_{i,j} = \int_{U_d} (z_i - \delta_{i,j})^+ \, \mu_j(\dd \bz) , \\
  \tBbeta = \int_{U_d} \bz \, \nu(\dd \bz) , \qquad
  \bC_k = \int_{U_d} \bz \bz^\top \mu_k(\dd \bz) , \qquad
  k \in \{1, \ldots, d\} ,
 \end{gather*}
 and \ $\obC$ \ is given by \eqref{obC} with the given
 \ $\bC_k$, \ $k \in \{1, \ldots, d\}$.
\end{Ex}

\begin{Ex} \
Let \ $(\bX_t)_{t\in\RR_+}$ \ be an irreducible and critical $d$-type
 ($d$-dimensional) CBI process with parameters
 \ $(d, \bc, \Bbeta, \bB, 0, \bzero)$.
\ Then
 \begin{gather*}
  \tbB = \bB , \qquad \tBbeta = \bbeta, \qquad
  \bC_k = 2 c_k \be_k \be_k^\top , \quad k \in \{1, \ldots, d\} , \qquad
  \obC = 2 \sum_{k=1}^d c_k (\be_k^\top \bu) (\be_k\be_k^\top).
 \end{gather*}
\end{Ex}

\section{Proof of Theorem \ref{main}}
\label{Proof_main}

The technique of the proof is somewhat similar to that of Theorem 3.1 in Isp\'any and Pap \cite{IspPap3}.
Let us introduce the notations
  \begin{gather*}
  \ttBbeta := \left( \int_0^1 \ee^{s\tbB} \, \dd s \right) \tBbeta
            \in \RR_+^d , \\
  \tbC := \sum_{k=1}^d
            \int_0^1
             (\be_k^\top \ee^{(1-s)\tbB} \bu) \,
             \ee^{s\tbB} \bC_k \ee^{s\tbB^\top}
             \! \dd s
       \in \RR_+^{d \times d}. 
  \end{gather*}
Note that \ $\ttBbeta = \EE(\bX_1)$ \ if \ $\bX_0 = \bzero$, \ see \eqref{EXbX}.

The process \ $(\bX_t)_{t\in\RR_+}$ \ is a time-homogeneous Markov process, hence
 Lemma 3.4 in Barczy et al. \cite{BarLiPap2} implies
 \[
   \EE(\bX_t \mid \bX_s = \bx)
   = \EE(\bX_{t-s} \mid \bX_0 = \bx)
   = \ee^{(t-s)\tbB} \bx
     + \left( \int_0^{t-s} \ee^{u\tbB} \, \dd u \right) \tBbeta
 \]
 for all \ $\bx \in \RR_+^d$ \ and \ $s, t \in \RR_+$ \ with \ $s < t$.
\ Using this formula, in order to prove \eqref{Conv_X}, let us introduce the
 sequence
 \begin{equation}\label{Mk}
  \begin{aligned}
   \bM_k &:= \bX_k - \EE(\bX_k \mid \cF^\bX_{k-1})
           = \bX_k - \EE(\bX_k \mid \bX_{k-1}) \\
         &\:= \bX_k -  \ee^{\tbB} \bX_{k-1}
              - \left( \int_0^1 \ee^{u\tbB} \, \dd u \right) \tBbeta
            = \bX_k -  \ee^{\tbB} \bX_{k-1} - \ttBbeta ,
   \qquad k \in \NN ,
  \end{aligned}
 \end{equation}
 which is a sequence of martingale differences with respect to the filtration
 \ $\bigl(\cF^\bX_k\bigr)_{k \in \ZZ_+}$, \ where, for each \ $t \in \RR_+$,
 \ $\cF^\bX_t$ \ denotes the \ $\sigma$-algebra generated by
 \ $(\bX_s)_{s\in[0,t]}$.
\ Consider the random step processes
 \begin{equation}\label{cMk}
  \bcM^{(n)}_t := n^{-1} \left( \bX_0 + \sum_{k=1}^{\nt} \bM_k \right) ,
  \qquad t \in \RR_+ , \qquad n \in \NN .
 \end{equation}
First we will verify convergence
 \begin{gather}\label{Conv_M}
  \bcM^{(n)} \distr \bcM \qquad \text{as \ $n \to \infty$,}
 \end{gather}
 where \ $(\bcM_t)_{t\in\RR_+}$ \ is the unique strong solution of the SDE
 \begin{equation}\label{SDE_M}
  \dd \bcM_t
  = \sqrt{\Big\langle \bv, \bcM_t + t \ttBbeta \Big\rangle^+ \tbC}
    \, \dd \bcW_t , \qquad t \in \RR_+ , \qquad \bcM_0 = \bzero ,
 \end{equation}
 where \ $(\bcW_t)_{t\in\RR_+}$ \ is a $d$-dimensional standard Brownian motion.
Here we note that the matrix
 \ $\Big\langle \bv, \bcM_t + t \ttBbeta \Big\rangle^+ \tbC$ \ is
 symmetric and positive semidefinite, and
 \ $\sqrt{\Big\langle \bv, \bcM_t + t \ttBbeta\Big\rangle^+ \tbC}$
 \ denotes its unique symmetric and positive semidefinite square root.
We want to apply Theorem \ref{Conv2DiffThm} with  \ $\bcU := \bcM$,
 \ $\bU^{(n)}_k := n^{-1} \bM_k$, \ $n, k \in \NN$,
 \ $\bU^{(n)}_0 := n^{-1} \bX_0$, \ $n \in \NN$, \ $\cF^{(n)}_k := \cF_k^\bX$
 \ for \ $n \in \NN$, \ $k \in \ZZ_+$, \ and with the coefficient function
 \ $\gamma : \RR_+ \times \RR^d \to \RR^{d \times d}$ \ of the SDE \eqref{SDE_M}
 given by
 \[
   \gamma(t, \bx)
   := \sqrt{\Big\langle\bv, \bx + t \ttBbeta \Big\rangle^+ \tbC} ,
   \qquad (t, \bx) \in \RR_+ \times \RR^d .
 \]
The aim of the following discussion is to show that the SDE \eqref{SDE_M} has
 a pathwise unique strong solution \ $\bigl(\bcM_t^{(\by_0)}\bigr)_{t\in\RR_+}$
 \ with initial value \ $\bcM_0^{(\by_0)} = \by_0$ \ for all \ $\by_0 \in \RR^d$.
\ First observe that if \ $\bigl(\bcM_t^{(\by_0)}\bigr)_{t\in\RR_+}$ \ is a strong
 solution of the SDE \eqref{SDE_M} with initial value
 \ $\bcM_0^{(\by_0)} = \by_0$, \ then, by It\^o's formula, the process
 \ $(\cP_t, \bcQ_t)_{t\in\RR_+}$, \ defined by
 \[
   \cP_t := \Big\langle \bv, \bcM_t^{(\by_0)}  + t \ttBbeta \Big\rangle , \qquad
   \bcQ_t := \bcM_t^{(\by_0)} - \cP_t \bu , \qquad t \in \RR_+ ,
 \]
 is a strong solution of the SDE
 \begin{equation}\label{SDE_d_Q}
  \begin{cases}
   \dd \cP_t
   = \bv^\top \ttBbeta \, \dd t
     + \sqrt{\cP_t^+} \, \bv^\top \sqrt{\tbC} \, \dd \bcW_t, \\[2mm]
   \dd \bcQ_t
   = - \bPi \ttBbeta \, \dd t
     + \sqrt{\cP_t^+} \, (\bI_d - \bPi) \sqrt{\tbC} \, \dd \bcW_t
  \end{cases}
  \qquad t \in \RR_+ ,
 \end{equation}
 with initial value
 \ $(\cP_0, \bcQ_0) = \big(\bv^\top \by_0, \, (\bI_d - \bPi)\by_0\bigr)$,
 \ where \ $\bPi := \bPi_{\tbB}$, \ see (iii) of Lemma \ref{FP}.
Indeed, the first SDE of \eqref{SDE_d_Q} is an easy consequence of the SDE
 \eqref{SDE_M}.
The second one can be checked as follows.
By It\^o's formula,
 \begin{align*}
  \dd \bcQ_t
  &= \dd \bcM_t^{(\by_0)} -  \bu \, \dd \cP_t
   = \dd \bcM_t^{(\by_0)}
     -  \bu \bv^\top (\dd \bcM_t^{(\by_0)} + \ttBbeta \, \dd t)
   = - \bPi \ttBbeta \, \dd t
     + (\bI_d - \bPi) \, \dd \bcM_t^{(\by_0)} \\
  &= - \bPi \ttBbeta \, \dd t
     + (\bI_d - \bPi)
       \sqrt{\Big\langle \bv, \bcM_t^{(\by_0)}  + t \ttBbeta \Big\rangle^+ \tbC}
       \, \dd \bcW_t \\
  &= - \bPi \ttBbeta \,\dd t
     + \sqrt{\cP_t^+} (\bI_d - \bPi) \sqrt{\tbC} \, \dd \bcW_t ,
   \qquad t \in \RR_+ ,
 \end{align*}
 where
 \ $\bcQ_0 = \by_0 - (\bv^\top \by_0) \bu = \by_0 - \bu \bv^\top \by_0
    = (\bI_d - \bPi) \by_0$.
\ Conversely, if
 \ $(\cP_t^{(p_0, \bq_0)}, \, \bcQ_t^{(p_0, \bq_0)})_{t\in\RR_+}$ \ is a strong
 solution of the SDE \eqref{SDE_d_Q} with initial value
 \ $\bigl(\cP_0^{(p_0, \bq_0)}, \, \bcQ_0^{(p_0, \bq_0)}\bigr) = (p_0, \bq_0)
    \in \RR \times \RR^d$,
 \ then, again by It\^o's formula,
 \[
   \bcM_t := \cP_t^{(p_0, \bq_0)} \, \bu + \bcQ_t^{(p_0, \bq_0)} ,
   \qquad t \in \RR_+ ,
 \]
 is a strong solution of the SDE \eqref{SDE_M} with initial value
 \ $\bcM_0 = p_0 \bu + \bq_0$.
\ The correspondence
 \ $\by_0 \leftrightarrow (p_0, \bq_0)
    := (\bv^\top \by_0, \, (\bI_d - \bPi)\by_0)$
 \ is a bijection between \ $\RR^d$ \ and
 \ $\RR \times \{\bq \in \RR^d : \bv^\top \bq = 0\}$, \ since
 \ $\by_0 = p_0 \bu + \bq_0$, \ and for all
 \ $(p_0, \bq_0) \in \RR \times \{\bq \in \RR^d : \bv^\top \bq = 0\}$, \ by
 (iii) of Lemma \ref{FP}, \ $\bv^\top(p_0 \bu + \bq_0) = p_0$, \ and
 \begin{align*}
  (\bI_d - \bPi)(p_0 \bu + \bq_0)
  = p_0 \bu + \bq_0 - p_0 \bPi \bu - \bPi \bq_0
  = p_0 \bu + \bq_0 - p_0 \bu \bv^\top \bu - \bu \bv^\top \bq_0
  = \bq_0 .
 \end{align*}
Hence it is enough to show that the SDE \eqref{SDE_d_Q} has a pathwise unique
 strong solution \ $(\cP_t^{(p_0, \bq_0)}, \, \bcQ_t^{(p_0, \bq_0)})_{t\in\RR_+}$
 \ with initial value
 \ $\bigl(\cP_0^{(p_0, \bq_0)}, \, \bcQ_0^{(p_0, \bq_0)}\bigr) = (p_0, \bq_0)$ \ for
 all \ $(p_0, \bq_0) \in \RR \times \{\bq \in \RR^d : \bv^\top \bq = 0\}$
 \ (actually, it turns out that it has a pathwise unique strong solution in
 case of any \ $(p_0, \bq_0) \in \RR \times \RR^d$).
\ The first equation of \eqref{SDE_d_Q} can be written in the form
 \begin{equation}\label{CIR3}
  \dd \cP_t = \bv^\top \ttBbeta \, \dd t
              + \sqrt{\cP_t^+} \, \sqrt{\langle \tbC \bv, \bv \rangle} \, \dd \cW_t ,
  \qquad t \in \RR_+ ,
 \end{equation}
 where \ $(\cW_t)_{t \in \RR_+}$ \ is a 1-dimensional standard Brownian motion.
Indeed,
 \ $\langle \tbC \bv, \bv \rangle = (\bv^\top \sqrt{\tbC}) (\bv^\top \sqrt{\tbC})^\top
    = \|\bv^\top \sqrt{\tbC}\|^2$,
 \ and hence if \ $\bv^\top \sqrt{\tbC} = \bzero \in \RR^{1\times d}$, \ then the
 above mentioned rewriting of the SDE \eqref{SDE_d_Q} is trivial, and if
 \ $\bv^\top \sqrt{\tbC} \ne \bzero \in \RR^{1\times d}$, \ then
 \ $\cW_t := \langle \tbC \bv, \bv \rangle^{-1/2} \bv^\top \sqrt{\tbC} \, \bcW_t$,
 \ $t \in \RR_+$, \ is a 1-dimensional standard Brownian motion.
Hence, by Remark \ref{REMARK1}, the first equation of the SDE \eqref{SDE_d_Q}
 has a pathwise unique strong solution \ $(\cP_t^{(p_0)})_{t\in\RR_+}$ \ with
 initial value \ $\cP_0^{(p_0)} = p_0$ \ for all \ $p_0 \in \RR$.
\ Clearly, the second equation of the SDE \eqref{SDE_d_Q} has a pathwise
 unique strong solution
 \[
   \bcQ_t^{(p_0, \bq_0)}
   = \bq_0 - \bPi \ttBbeta \, t
     + (\bI_d - \bPi)
       \sqrt{\tbC} \int_0^t \sqrt{(\cP_s^{(p_0)})^+} \, \dd \bcW_s ,
   \qquad t \in \RR_+ ,
 \]
 with initial value \ $\bcQ_0^{(p_0, \bq_0)} = \bq_0$ \ for all
 \ $(p_0, \bq_0) \in \RR \times \RR^d$.
\ Consequently, the SDE \eqref{SDE_d_Q}, and hence the SDE \eqref{SDE_M} admit
 a pathwise unique strong solution with an arbitrary initial value.

Now we show that conditions (i) and (ii) of Theorem \ref{Conv2DiffThm} hold.
We have to check that for each \ $T > 0$,
 \begin{align} \label{Cond1}
  &\sup_{t\in[0,T]}
    \bigg\| \frac{1}{n^2} \sum_{k=1}^{\nt} \var(\bM_k \mid \cF^\bX_{k-1})
            - \left( \int_0^t (\bcR^{(n)}_s)^+ \, \dd s \right)
              \tbC \bigg\|
   \stoch 0,\\
  &\frac{1}{n^2}
   \sum_{k=1}^{\nT}
    \EE(\|\bM_k\|^2 \bbone_{\{\|\bM_k\| > n \theta\}} \mid \cF^\bX_{k-1})
   \stoch 0   \qquad\text{for all \ $\theta>0$} \label{Cond2}
 \end{align}
 as \ $n \to \infty$, \ where the process \ $(\bcR^{(n)}_s)_{s\in\RR_+}$ \ is
 defined by
 \begin{equation}\label{Rnt}
  \bcR^{(n)}_s := \bv^\top \Bigl(\bcM^{(n)}_s + s \ttBbeta\Bigr) , \qquad
  s \in \RR_+ , \qquad n \in \NN .
 \end{equation}
By \eqref{Mk},
 \begin{align*}
  \bcR^{(n)}_s
  & = \bv^\top
      \bigg( n^{-1}\bX_0 + n^{-1}
             \sum_{k=1}^\ns (\bX_k - \ee^{\tbB} \bX_{k-1} - \ttBbeta)
             + s \ttBbeta \bigg) \\
  & = n^{-1}\bv^\top\bX_0 + n^{-1} \sum_{k=1}^\ns
             \big( \bv^\top \bX_k - \bv^\top \ee^{\tbB} \bX_{k-1}
                   - \bv^\top \ttBbeta \big)
             + s \bv^\top \ttBbeta \\
  & = n^{-1}\bv^\top\bX_0 + n^{-1} \sum_{k=1}^\ns
             \big( \bv^\top \bX_k - \bv^\top \bX_{k-1} - \bv^\top \ttBbeta \big)
             + s \bv^\top \ttBbeta \\
  & = n^{-1}\bv^\top \bX_\ns
      + \left( s - \frac{\ns}{n} \right) \bv^\top \ttBbeta ,
 \end{align*}
 where we used that \ $\bv$ \ is a left eigenvector of \ $\ee^{\tbB}$
 \ belonging to the eigenvalue 1 (due to the fact that \ $\bv$ \ is a left
 eigenvector of \ $\tbB$ \ corresponding to the eigenvalue \ $s(\tbB) = 0$).
\ Thus \ $(\bcR^{(n)}_s)^+ = \bcR^{(n)}_s$, \ and
 \begin{align*}
  \int_0^t (\bcR^{(n)}_s)^+ \, \dd s
  &= \frac{1}{n^2} \sum_{\ell=0}^{\nt-1} \bv^\top \bX_\ell
     + \left(t - \frac{\nt}{n}\right) \frac{1}{n} \bv^\top \bX_{\nt} \\
  &\quad
     + \frac{t^2}{2}\bv^\top \ttBbeta
     - \frac{1}{n^2} \Biggl(\sum_{\ell=1}^{\nt-1}\ell \Biggr) \bv^\top \ttBbeta
     - \left(t - \frac{\nt}{n}\right)\frac{\nt}{n} \bv^\top \ttBbeta \\
  &= \frac{1}{n^2} \sum_{\ell=0}^{\nt-1} \bv^\top \bX_\ell
     + \frac{nt-\nt}{n^2} \bv^\top \bX_{\nt}
     + \frac{\nt+(nt-\nt)^2}{2n^2} \bv^\top \ttBbeta .
 \end{align*}
We have \ $\var(\bM_k \mid \cF^\bX_{k-1}) = \var(\bX_k \mid \bX_{k-1})$ \ and
 \ $\var(\bX_k \mid \bX_{k-1} = \bx) = \var(\bX_1 \mid \bX_0 = \bx)$ \ for all
 \ $\bx \in \RR_+^d$, \ since \ $(\bX_t)_{t\in\RR_+}$ \ is a time-homogeneous
 Markov process.
Hence Lemma 4.4 in Barczy et al. \cite{BarLiPap3} implies
 \[
   \frac{1}{n^2} \sum_{k=1}^{\nt} \var(\bM_k \mid \cF^\bX_{k-1})
   = \frac{\nt}{n^2} \bV
     + \frac{1}{n^2}
       \sum_{k=1}^{\nt}
        \sum_{\ell=1}^d
         \int_0^1
          (\be_\ell^\top \ee^{(1-s)\tbB} \bX_{k-1}) \,
          \ee^{s\tbB} \bC_\ell \ee^{s\tbB^\top} \!
          \dd s ,
 \]
 where
 \begin{align*}
  \bV := \int_0^1
           \ee^{u\tbB}
           \left( \int_{U_d} \bz \bz^\top \nu(\dd \bz) \right)
           \ee^{u\tbB^\top} \!
           \dd u
          + \sum_{k=1}^d
             \int_0^1
              \left( \int_0^{1-u} \be_k^\top \ee^{v\tbB} \tBbeta \, \dd v \right)
              \ee^{u\tbB} \bC_k \ee^{u\tbB^\top} \! \dd u .
 \end{align*}
Hence, in order to show \eqref{Cond1}, it suffices to prove
 \begin{equation}\label{Cond11}
  \begin{split}
   &n^{-2}
    \sup_{t \in [0,T]}
     \sum_{k=0}^{\nt-1}
      \left\| \sum_{\ell=1}^d
               \int_0^1
                (\be_\ell^\top \ee^{(1-s)\tbB} \bX_k) \,
                \ee^{s\tbB} \bC_\ell \ee^{s\tbB^\top} \!
                \dd s
              - (\bv^\top \bX_k) \tbC \right\|
    \stoch 0 , \\
   &n^{-2} \sup_{t \in [0,T]} \|\bX_{\nt}\| \stoch 0
  \end{split}
 \end{equation}
 as \ $n \to \infty$.
\ Observe that
 \begin{align*}
   (\bv^\top \bX_k) \tbC
  & = \sum_{\ell=1}^d
      \int_0^1
       (\be_\ell^\top \ee^{(1-s)\tbB} \bu) ( \bv^\top \bX_k) \,
       \ee^{s\tbB} \bC_\ell \ee^{s\tbB^\top} \!
       \dd s \\
  & = \sum_{\ell=1}^d
      \int_0^1
       (\be_\ell^\top \ee^{(1-s)\tbB} \bPi \bX_k) \,
       \ee^{s\tbB} \bC_\ell \ee^{s\tbB^\top} \!
       \dd s , \qquad k \in \ZZ_+ .
 \end{align*}
From \eqref{Mk} we obtain the recursion
 \begin{equation*}
  \bX_k = \ee^{\tbB} \bX_{k-1} + \bM_k + \ttBbeta , \qquad k \in \NN ,
 \end{equation*}
 hence we conclude
 \begin{equation}\label{X}
  \bX_k = \ee^{k\tbB} \bX_0
          + \sum_{j=1}^k \ee^{(k-j)\tbB} \Bigl(\bM_j + \ttBbeta\Bigr) ,
  \qquad k \in \NN .
 \end{equation}
Using \eqref{X}, for all \ $k \in \ZZ_+$, \ we obtain
 \begin{multline*}
  \sum_{\ell=1}^d
   \int_0^1
    (\be_\ell^\top \ee^{(1-s)\tbB} \bX_k) \,
    \ee^{s\tbB} \bC_\ell \ee^{s\tbB^\top} \!
    \dd s
  - (\bv^\top \bX_k) \tbC \\
  \begin{aligned}
   &= \sum_{\ell=1}^d
       \int_0^1
        \left[ \be_\ell^\top \ee^{(1-s)\tbB} (\bI_d - \bPi) \bX_k \right]
        \ee^{s\tbB} \bC_\ell \ee^{s\tbB^\top} \!
        \dd s \\
   &= \sum_{\ell=1}^d
       \int_0^1
        \left[ \be_\ell^\top \ee^{(1-s)\tbB} (\bI_d - \bPi)
               \left\{ \ee^{k\tbB} \bX_0
                       + \sum_{j=1}^k
                          \ee^{(k-j)\tbB} (\bM_j + \ttBbeta) \right\} \right]
        \ee^{s\tbB} \bC_\ell \ee^{s\tbB^\top} \!
        \dd s  \\
   &=  \sum_{\ell=1}^d
       \int_0^1
        \left[ \be_\ell^\top \ee^{(1-s)\tbB}
               \left\{ (\ee^{k\tbB} - \bPi) \bX_0
                       + \sum_{j=1}^k \left( \ee^{(k-j)\tbB} - \bPi \right)
                         (\bM_j + \ttBbeta) \right\} \right]
        \ee^{s\tbB} \bC_\ell \ee^{s\tbB^\top} \!
        \dd s ,
  \end{aligned}
 \end{multline*}
 since, by (iii) of Lemma \ref{FP},
 \ $\bPi \ee^{\tbB} = \left(\lim_{t\to\infty} \ee^{t\tbB}\right) \ee^{\tbB}
    = \lim_{t\to\infty} \ee^{(t+1)\tbB} = \bPi$
 \ implies \ $(\bI_d - \bPi) \ee^{(k-j)\tbB} = \ee^{(k-j)\tbB} - \bPi$.
\ Hence, by (iv) of Lemma \ref{FP},
 \begin{multline*}
  \sum_{k=0}^\nt
   \left\| \sum_{\ell=1}^d
            \int_0^1
            (\be_\ell^\top \ee^{(1-s)\tbB} \bX_k) \,
            \ee^{\tbB s} \bC_\ell \ee^{s\tbB^\top} \!
            \dd s
           - (\bv^\top \bX_k) \tbC \right\| \\
  \begin{aligned}
   & \leq c
          \sum_{k=0}^\nt
           \sum_{\ell=1}^d
             \left\{ \ee^{-k\kappa} \|\bX_0\|
                    + \sum_{j=1}^k
                       \ee^{-(k-j)\kappa} \| \bM_j + \ttBbeta \| \right\}
            \int_0^1
             \|\ee^{(1-s)\tbB}\| \|\ee^{s\tbB} \bC_\ell \ee^{s\tbB^\top} \! \|
             \, \dd s \\
   & \leq c (c + \|\bPi\|)^3
          \sum_{\ell=1}^d
           \|\bC_\ell\|
           \left\{ \|\bX_0\| \sum_{k=0}^{\lfloor nt\rfloor}\ee^{-k\kappa}
                   + \sum_{j=1}^\nt \sum_{k=j}^\nt
                      \ee^{-(k-j)\kappa} (\|\bM_j\| + \|\ttBbeta\|) \right\} \\
   & \leq \frac{c (c + \|\bPi\|)^3}{1 - \ee^{-\kappa}}
          \bigg( \|\bX_0\| +  \nt \cdot \|\ttBbeta\|
                 + \sum_{j=1}^\nt \|\bM_j\| \bigg)
          \sum_{\ell=1}^d \|\bC_\ell\| ,
  \end{aligned}
 \end{multline*}
 since
 \ $\sum_{k=0}^{\lfloor nt\rfloor} \ee^{-k\kappa} < \sum_{k=0}^\infty \ee^{-k\kappa}
    = \frac{1}{1-\ee^{-\kappa}}$,
 \ and, by (iv) of Lemma \ref{FP},
 \begin{equation}\label{c+Pi}
  \|\ee^{t\tbB^\top} \! \| = \|\ee^{t\tbB}\| \leq c + \|\bPi\| ,
  \qquad t \in \RR_+ ,
 \end{equation}
 where \ $c := c_{\tbB} \in \RR_{++}$ \ and \ $\kappa := \kappa_{\tbB} \in \RR_{++}$ \ are given in (iv) of Lemma \ref{FP}.

Moreover, by \eqref{X} and \eqref{c+Pi},
 \begin{align*}
   \|\bX_{\nt}\|
  & \leq  \| \ee^{\nt\tbB} \| \|\bX_0\|
        + \sum_{j=1}^{\nt}
           \| \ee^{(\nt-j)\tbB} \| \|\bM_j + \ttBbeta\|\\
  & \leq (c + \|\bPi\|)
        \bigg( \|\bX_0\| +  \nt \cdot \|\ttBbeta\|
               + \sum_{j=1}^{\nt} \|\bM_j\| \bigg) .
 \end{align*}
Consequently, in order to prove \eqref{Cond11}, it suffices to show
 \[
  \frac{1}{n^2} \sum_{j=1}^{\nT} \|\bM_j\| \stoch 0 \qquad
  \text{as \ $n \to \infty$.}
 \]
In fact, Lemma \ref{moment_estimations_1_2} yields
 \ $n^{-2} \sum_{j=1}^{\nT} \EE(\|\bM_j\|) \leq n^{-2} \sum_{j=1}^{\nT} \sqrt{\EE(\|\bM_j\|^2)} \to 0$
 \ as \ $n \to \infty$, \ thus we obtain \eqref{Cond1}.

Next we check condition \eqref{Cond2}.
We have
 \[
   \EE(\|\bM_k\|^2 \bbone_{\{\|\bM_k\| > n \theta\}} \mid \cF^\bX_{k-1})
   \leq n^{-2} \theta^{-2} \EE(\|\bM_k\|^4 \mid \cF^\bX_{k-1}) .
 \]
Moreover, \ $n^{-4} \sum_{k=1}^{\nT} \EE\bigl(\|\bM_k\|^4\bigr) \to 0$ \ as
 \ $n \to \infty$, \ since \ $\EE(\|\bM_k\|^4) = \OO(k^2)$ \ by Lemma
 \ref{moment_estimations_1_2} (at this point we used the moment condition
 \ $\EE(\|\bX_0\|^4) < \infty$).
Thus we obtain \eqref{Cond2}, and hence, convergence \eqref{Conv_M}.

Applying a version of the continuous mapping theorem together with
 \eqref{Conv_M} and \eqref{X}, we will show
 \begin{gather}\label{Conv_bX}
  (\bcX^{(n)}_t)_{t\in\RR_+} \distr (\bcX_t)_{t\in\RR_+}
  \qquad \text{as \ $n \to \infty$,}
 \end{gather}
 where
 \begin{align}\label{help1}
   \bcX_t := \bPi \Bigl(\bcM_t + t \ttBbeta\Bigr) , \qquad t \in \RR_+ ,
 \end{align}
  where \ $(\bcM_t)_{t\in\RR_+}$ \ is defined in \eqref{SDE_M}.
We want to apply Lemma \ref{Conv2Funct}.
By \eqref{X}, \ $\bcX^{(n)} = \Psi_n(\bcM^{(n)})$, $n\in\NN$, \ where the
 mapping \ $\Psi_n : \DD(\RR_+, \RR^d) \to \DD(\RR_+, \RR^d)$ \ is given by
 \[
   \Psi_n(f)(t)
   := \ee^{\nt \tbB} f(0)
      + \sum_{j=1}^\nt
         \ee^{(\nt-j)\tbB}
         \left( f\left(\frac{j}{n}\right) - f\left(\frac{j-1}{n}\right)
                + n^{-1} \ttBbeta \right)
 \]
 for \ $f \in \DD(\RR_+, \RR^d)$, \ $t \in \RR_+$, \ $n \in \NN$.
\ Further, by \eqref{help1}, \ $\bcX = \Psi(\bcM)$, \ where the mapping
 \ $\Psi : \DD(\RR_+, \RR^d) \to \DD(\RR_+, \RR^d)$ \ is given by
 \[
   \Psi(f)(t) := \bPi \Bigl(f(t) + t \ttBbeta\Bigr) , \qquad
   f \in \DD(\RR_+, \RR^d), \qquad t \in \RR_+ .
 \]
Measurability of the mappings \ $\Psi_n$, \ $n \in \NN$, \ and \ $\Psi$ \ can
 be checked similarly as in Barczy et al. \cite[page 603]{BarIspPap0}.
We only note that, with the notations of Barczy et al. \cite{BarIspPap0}, for
 all \ $n,N\in\NN$, \ the mappings
 \ $\psi_n^{N,1} : \DD(\RR_+,\RR^d) \to \RR^{(nN+1)d}$ \ and
 \ $\psi_n^{N,2} : \RR^{(nN+1)d} \to \DD(\RR_+,\RR^d)$ \ should be defined by
 \begin{align*}
  &\psi_n^{N,1}(f)
   := \left(f(0), f\left(\frac{1}{n}\right), f\left(\frac{2}{n}\right),
            \ldots, f(N)\right) , \\
  &\psi_n^{N,2}(\bx_0, \bx_1, \ldots, \bx_{nN})(t)
   := \ee^{\lfloor nt\rfloor\tbB} \bx_0
      + \sum_{j=1}^{\lfloor nt\rfloor}
         \ee^{(\lfloor nt\rfloor - j)\tbB}
         (\bx_j - \bx_{j-1} + n^{-1} \ttBbeta)
 \end{align*}
 for \ $f \in \DD(\RR_+,\RR^d)$, \ $t \in \RR_+$ \ and
 \ $(\bx_0^\top, \bx_1^\top, \ldots, \bx_{nN}^\top)^\top \in \RR^{(nN+1)d}$.

The aim of the following discussion is to show that the set
 \ $C := \{ f \in \CC(\RR_+, \RR^d) : \bPi f(0) = f(0) \}$ \ satisfies
 \ $C \in \cD_\infty(\RR_+, \RR^d)$, \ $C \subset C_{\Psi, \, (\Psi_n)_{n \in \NN}}$
 \ and \ $\PP(\bcM \in C) = 1$.

First note that
 \ $C = \CC(\RR_+, \RR^d)
        \cap \pi_0^{-1}\bigl( \{ \bx \in \RR^d
                                : (\bI_d - \bPi) \bx = \bzero \} \bigr)$,
 \ where \ $\pi_0 : \DD(\RR_+, \RR^d) \to \RR^d$ \ denotes the projection
 defined by \ $\pi_0(f) := f(0)$ \ for \ $f \in \DD(\RR_+, \RR^d)$.
\ Using that \ $\CC(\RR_+, \RR^d) \in \cD_\infty(\RR_+, \RR^d)$ \ (see, e.g.,
 Ethier and Kurtz \cite[Problem 3.11.25]{EthKur}), the mapping
 \ $\RR^d \ni \bx \mapsto (\bI_d - \bPi) \bx \in \RR^d$ \ is measurable and
 that \ $\pi_0$ \ is measurable (see, e.g., Ethier and Kurtz
 \cite[Proposition 3.7.1]{EthKur}), we obtain
 \ $C \in \cD_\infty(\RR_+, \RR^d)$.

Fix a function \ $f \in C$ \ and a sequence \ $(f_n)_{n\in\NN}$ \ in
 \ $\DD(\RR_+, \RR^d)$ \ with \ $f_n \lu f$ \ as \ $n\to\infty$.
\ By the definition of \ $\Psi$, \ we have \ $\Psi(f) \in \CC(\RR_+, \RR^d)$.
\ Furthermore,
 \begin{align*}
  \Psi_n(f_n)(t)
  &= \bPi
     \left( f_n\left(\frac{\nt}{n}\right)
            + \frac{\nt}{n} \ttBbeta \right)
     + (\ee^{\nt \tbB} - \bPi) f_n(0) \\
  &\quad
     + \sum_{j=1}^{\nt}
        \bigl(\ee^{(\nt-j)\tbB} - \bPi\bigr)
        \left( f_n\left(\frac{j}{n}\right)
               - f_n\left(\frac{j-1}{n}\right)
               + \frac{1}{n} \ttBbeta \right) ,
 \end{align*}
 hence we have for all \ $t\in\RR_+$,
 \begin{align*}
  \|\Psi_n(f_n)(t)-\Psi(f)(t)\|
  &\leq \| \bPi \|
        \left( \left\| f_n\left(\frac{\nt}{n}\right)
                       - f(t) \right\|
               + \frac{1}{n} \|\ttBbeta\| \right)
        + \|(\ee^{\nt \tbB} - \bPi) f_n(0)\| \\
  &\quad
        + \sum_{j=1}^{\nt}
           \bigl\| \ee^{(\nt-j)\tbB} - \bPi \bigr\|
           \left( \left\| f_n\left(\frac{j}{n}\right)
                          - f_n\left(\frac{j-1}{n}\right) \right\|
                  + \frac{1}{n} \|\ttBbeta\| \right) .
 \end{align*}
Here for all \ $T > 0$ \ and \ $t \in [0, T]$,
 \begin{align*}
  \left\| f_n\left(\frac{\nt}{n}\right) - f(t) \right \|
  &\leq \left\| f_n\left(\frac{\nt}{n}\right)
                - f\left(\frac{\nt}{n}\right) \right\|
        + \left\| f\left(\frac{\nt}{n}\right) - f(t) \right\| \\
  &\leq \sup_{t \in [0,T]} \|f_n(t) - f(t)\| + \omega_T(f, n^{-1}) ,
 \end{align*}
 where \ $\omega_T(f, \cdot)$ \ is the modulus of continuity of \ $f$ \ on
 \ $[0, T]$, \ and we have \ $\omega_T(f, n^{-1}) \to 0$ \ since \ $f$ \ is
 continuous (see, e.g., Jacod and Shiryaev \cite[VI.1.6]{JacShi}).
In a similar way, for all \ $T > 0$,
 \[
   \left\| f_n\left(\frac{j}{n}\right)
           - f_n\left(\frac{j-1}{n}\right) \right\|
   \leq \omega_T(f, n^{-1}) + 2 \sup_{t \in [0,T]} \|f_n(t) - f(t)\| ,
   \qquad  j \in \{1, \ldots, \nT\} .
 \]
By (iv) of Lemma \ref{FP}, for all \ $T>0$ \ and \ $t\in[0,T]$,
 \[
   \sum_{j=1}^{\nt} \bigl\| \ee^{(\nt-j)\tbB} - \bPi \bigr\|
   \leq \sum_{j=1}^\nt c \ee^{-(\nt-j)\kappa}
   \leq \frac{c}{1 - \ee^{-\kappa}} .
 \]
Further, for all \ $T > 0$ \ and \ $t \in [0, T]$, \ by (iv) of Lemma \ref{FP},
 \begin{align*}
  \|(\ee^{\nt \tbB} - \bPi) f_n(0)\|
  &\leq \|(\ee^{\nt \tbB} - \bPi) (f_n(0) - f(0))\|
        + \|(\ee^{\nt \tbB} - \bPi) f(0)\| \\
  &\leq c \sup_{t\in[0,T]} \|f_n(t) - f(t)\| ,
 \end{align*}
 since \ $f \in C$ \ implies \ $(\ee^{\nt \tbB} - \bPi) f(0) = \bzero$.
\ Indeed, by part (iii) of Lemma \ref{FP},
 \ $\ee^{\nt \tbB} \bPi = \ee^{\nt \tbB} \lim_{s \to \infty} \ee^{s \tbB}
    = \lim_{s \to \infty} \ee^{(\nt+s) \tbB} = \bPi$,
 \ hence \ $\ee^{\nt \tbB} f(0) = \ee^{\nt \tbB} \bPi f(0) = \bPi f(0)$.
\ Using that \ $f_n \lu f$ \ as \ $n \to \infty$, \ we have
 \ $\Psi_n(f_n) \lu \Psi(f)$ \ as \ $n \to \infty$.
\ Thus we conclude \ $C \subset C_{\Psi, \, (\Psi_n)_{n \in \NN}}$.

By the definition of a weak solution (see, e.g., Karatzas and Shreve
 \cite[Definition 3.1, Section 5.3]{KarShr}), \ $\bcM$ \ has continuous sample
 paths almost surely, hence, since \ $\bcM_0=\bzero$, \ we have
 \ $\PP(\bcM \in C) = 1$.
\ Consequently, by Lemma \ref{Conv2Funct}, we obtain
 \ $\bcX^{(n)} = \Psi_n(\bcM^{(n)}) \distr \Psi(\bcM) = \bcX$ \ as
 \ $n \to \infty$.

It remains to show that the limit process given by \eqref{help1} coincides in
 law with the corresponding one in \eqref{Conv_X}.
Using \ $\bPi = \bu \bv^\top$ \ and \ $\bv^\top \bu = 1$, \ we get that the
 process \ $\cX_t := \bv^\top \bcX_t$, \ $t \in \RR_+$ \ (where
 \ $(\bcX_t)_{t\in\RR_+}$ \ is given by \eqref{help1}) satisfies
 \[
   \cX_t
   = \bv^\top \bPi \Bigl(\bcM_t + t \ttBbeta\Bigr)
   = \bv^\top \Bigl(\bcM_t + t \ttBbeta\Bigr) ,
   \qquad t \in \RR_+ ,
 \]
 hence
 \ $\cX_t \bu = \bu \bv^\top \Bigl(\bcM_t + t \ttBbeta\Bigr)
    = \bPi \Bigl(\bcM_t + t \ttBbeta\Bigr) = \bcX_t$,
 \ $t \in \RR_+$.
\ By \eqref{SDE_M} and It\^o's formula we obtain that \ $(\cX_t)_{t\in\RR_+}$
 \ is a strong solution of the SDE
 \[
   \dd \cX_t
   = \bv^\top \ttBbeta \, \dd t
     + \sqrt{\cX_t^+} \bv^\top \sqrt{\tbC} \, \dd \bcW_t , \qquad t \in \RR_+ ,
   \qquad \cX_0 = 0 ,
 \]
 where \ $(\bcW_t)_{t\in\RR_+}$ \ is a $d$-dimensional standard Brownian motion.
This equation can be written in the form
 \begin{align}\label{SDE_hullamos}
  \dd\cX_t = \bv^\top \ttBbeta \, \dd t + \sqrt{\bv^\top \tbC \bv\cX_t^+}\,\dd\cW_t,
  \qquad t\in\RR_+,\qquad \cX_0=0,
 \end{align}
 with some 1-dimensional standard Brownian motion \ $(\cW_t)_{t\in\RR_+}$.
\ Indeed,
 \ $\bv^\top \tbC \bv = (\bv^\top \sqrt{\tbC}) (\bv^\top \sqrt{\tbC})^\top
    = \|\bv^\top \sqrt{\tbC}\|^2$,
 \ and hence if \ $\bv^\top \sqrt{\tbC} = \bzero\in\RR^{1\times d}$, \ then the
 above mentioned rewriting of the SDE in question is trivial, and if
 \ $\bv^\top \sqrt{\tbC} \ne \bzero \in \RR^{1\times d}$, \ then
 \ $\cW_t := (\bv^\top \tbC \bv)^{-1/2} \bv^\top \sqrt{\tbC} \, \bcW_t$,
 \ $t \in \RR_+$, \ is a 1-dimensional standard Brownian motion.
Finally, the SDE \eqref{SDE_hullamos} can be written in the form \eqref{SDE_X},
 since
 \begin{align*}
  \bv^\top \ttBbeta
    = \int_0^1 \bv^\top\ee^{s\tbB}\tBbeta\,\dd s
    = \int_0^1 \bv^\top \tBbeta\,\dd s
    = \bv^\top \tBbeta,
 \end{align*}
 and
 \begin{align*}
  \langle \tbC\bv,\bv\rangle
  = \bv^\top \tbC\bv
  & = \sum_{k=1}^d (\be_k^\top \bu) \int_0^1 \bv^\top \ee^{s\tbB}\bC_k \ee^{s\tbB^\top}\bv\,\dd s
    = \sum_{k=1}^d (\be_k^\top \bu) \int_0^1 \bv^\top \bC_k \bv\,\dd s \\
  & = \sum_{k=1}^d (\be_k^\top\bu) \bv^\top \bC_k \bv
    = \bv^\top \left( \sum_{k=1}^d  (\be_k^\top\bu) \bC_k \right) \bv
    = \bv^\top \obC \bv
    = \langle \obC \bv, \bv \rangle .
 \end{align*}
Hence \ $(\cX_t)_{t\in\RR_+}$ \ is a strong solution of the SDE \eqref{SDE_X} and
 consequently, we conclude \eqref{Conv_X}.
\proofend

\vspace*{5mm}

\appendix

\vspace*{5mm}

\noindent{\bf\Large Appendices}

\section{Frobenius--Perron type results}
\label{section_moments}

For the classification of CBI processes and for the estimation of the moments
 of a CBI process \ $(\bX_t)_{t\in\RR_+}$, \ we need some Frobenius--Perron type
 statements about the asymptotic behaviour of \ $\ee^{t\bA}$ \ as
 \ $t \to \infty$, \ where
 \ $\bA = (a_{i,j})_{i,j\in\{1,\ldots,d\}} \in \RR^{d \times d}_{(+)}$, \ i.e., \ $\bA$
 \ is essentially non-negative.
Note that then \ $\ee^\bA \in \RR_+^{d \times d}$.
\ Indeed, we have \ $\bA - a_\bA \bI_d \in \RR_+^{d \times d}$ \ with
 \ $a_\bA := \min_{i\in\{1,\ldots,d\}} a_{i,i}$, \ thus
 \ $\ee^\bA = \ee^{a_\bA} \ee^{\bA - a_\bA \bI_d} \in \RR_+^{d \times d}$.
\ Recall that \ $s(\bA) = \max_{\lambda\in\sigma(\bA)} \Re(\lambda)$, \ where
 \ $\sigma(\bA)$ \ denotes the spectrum of \ $\bA$, \ i.e., the set of the
 eigenvalues of \ $\bA$.

\begin{Lem}\label{FP1}
Suppose that \ $\bA \in \RR^{d \times d}_{(+)}$.
\ Then the following statements are equivalent:
\renewcommand{\labelenumi}{{\rm(\roman{enumi})}}
\begin{enumerate}
 \item
  there exists \ $t_0 \in \RR_{++}$ \ such that
   \ $\ee^{t_0\bA} \in \RR^{d \times d}_{++}$;
 \item
  for all \ $t \in \RR_{++}$, \ we have \ $\ee^{t\bA} \in \RR^{d \times d}_{++}$;
 \item
  $\bA$ \ is irreducible.
\end{enumerate}
\end{Lem}

\noindent
\textbf{Proof.}
If \ $d = 1$, \ then the statement is trivial.
If \ $d \geq 2$, \ then the statement follows by Berman and Plemmons
 \cite[Chapter 6, Theorem 3.12]{BerPle}.
However, for the claim (i) $\Longrightarrow$ (ii) we give an independent
 proof, which may be interesting on its own.
Let
 \ $\tbA = (\ta_{i,j})_{i,j\in\{1,\ldots,d\}}
    := \bA - a_\bA \bI_d \in \RR_+^{d \times d}$
 \ and \ $K := 1 + \max_{i\in\{1,\ldots,d\}} \be_i^\top \tbA \bone \in \RR_{++}$,
 \ where \ $\bone := (1, \ldots, 1)^\top \in \RR^d$.
\ Then the matrix
 \ $\bQ = (q_{i,j})_{i,j\in\{1,\ldots,d+1\}} \in \RR_+^{(d+1)\times(d+1)}$, \ given by
 \[
   q_{i,j}
   := \begin{cases}
       K^{-1} \, \ta_{i,j} , & \text{if \ $i, j \in \{1, \ldots, d\}$,} \\
       1 - K^{-1} \be_i^\top \tbA \bone ,
        & \text{if \ $i \in \{1, \ldots, d\}$ \ and \ $j = d+1$,} \\
       0 , & \text{if \ $i = d+1$ \ and \ $j \in \{1, \ldots, d\}$,} \\
       1 , & \text{if \ $i = d+1$ \ and \ $j = d+1$,}
      \end{cases}
 \]
 is a stochastic matrix, since the entries are non-negative and
 \ $\tbe_i^\top \bQ \tbone = 1$ \ for all
 \ $i \in \{1, \ldots, d+1\}$, \ where \ $\tbe_1$, \ldots, $\tbe_{d+1}$
 \ denotes the natural basis in \ $\RR^{d+1}$, \ and
 \ $\tbone := (1, \ldots, 1)^\top \in \RR^{d+1}$.
\ Indeed, for each \ $i \in \{1, \ldots, d\}$, \ we have
 \ $\tbe_i^\top \bQ \tbone
    = \be_i^\top K^{-1} \tbA \bone + (1 - K^{-1} \be_i^\top \tbA \bone) = 1$,
 \ and \ $\tbe_{d+1}^\top \bQ \tbone = q_{d+1,d+1} = 1$.
\ By Chung \cite[Theorem II.1.5]{Chu1}, for all \ $i, j \in \{1, \ldots, d+1\}$,
 \ there are two possibilities, namely, either
 \ $\tbe_i^\top \ee^{t\bQ} \, \tbe_j > 0$ \ for all \ $t \in \RR_{++}$, \ or
 \ $\tbe_i^\top \ee^{t\bQ} \, \tbe_j = 0$ \ for all \ $t \in \RR_{++}$.
\ Clearly, \ $\ta_{i,j} = K q_{i,j}$, \ $i, j \in \{1, \ldots, d\}$, \ and
 \ $\tbe_{d+1}^\top \bQ \, \tbe_j = 0$, \ $j \in \{1, \ldots, d\}$, \ imply
 \ $\be_i^\top \tbA^n \be_j = K^n \, \tbe_i^\top \bQ^n \, \tbe_j$ \ for all
 \ $i, j \in \{1, \ldots, d\}$ \ and \ $n \in \NN$.
\ Indeed, this obviously holds for \ $n = 1$, \ and, by induction,
 \begin{align*}
  \be_i^\top \tbA^{n+1} \be_j
  &= \sum_{k=1}^d (\be_i^\top \tbA^n \be_k) (\be_k^\top \tbA \be_j)
   = K^{n+1} \sum_{k=1}^d
             (\tbe_i^\top \bQ^n \, \tbe_k) (\tbe_k^\top \bQ \, \tbe_j) \\
  &= K^{n+1} \sum_{k=1}^{d+1}
             (\tbe_i^\top \bQ^n \, \tbe_k) (\tbe_k^\top \bQ \, \tbe_j)
   = K^{n+1} \tbe_i^\top \bQ^n \, \tbe_j .
 \end{align*}
Consequently,
 \[
   \be_i^\top \ee^{t\tbA} \be_j
   = \be_i^\top \sum_{n=0}^\infty \frac{t^n \tbA^n}{n!} \be_j
   = \tbe_i^\top \sum_{n=0}^\infty \frac{K^n t^n \bQ^n}{n!} \tbe_j
   = \tbe_i^\top \ee^{K t \bQ} \, \tbe_j
 \]
 for all \ $i, j \in \{1, \ldots, d\}$.
\ Now the assumption \ $\ee^{t_0\bA} \in \RR^{d \times d}_{++}$ \ implies
 \ $\tbe_i^\top \ee^{K t_0 \bQ} \, \tbe_j = \be_i^\top \ee^{t_0\tbA} \be_j
    = \ee^{-t_0a_\bA} \be_i^\top \ee^{t_0\bA} \be_j > 0$
 \ for all \ $i, j \in \{1, \ldots, d\}$, \ hence, by Chung
 \cite[Theorem II.1.5]{Chu1},
 \ $\be_i^\top \ee^{t\tbA} \be_j = \tbe_i^\top \ee^{K t \bQ} \, \tbe_j > 0$ \ for all
 \ $i, j \in \{1, \ldots, d\}$ \ and \ $t\in\RR_{++}$, \ and we conclude (ii).
\proofend

\begin{Rem}\label{FP1_Rem}
Exercise 7.7.4 of the Internet Seminar \cite{IS} claims that the above
 statements are equivalent to irreducibility of the matrix \ $\ee^{t\bA}$ \ for
 some or for any \ $t \in \RR_{++}$.
\ Thus a multi-type CBI process is irreducible if and only if \ $\ee^{\tbB}$
 \ is irreducible.
\proofend
\end{Rem}

\begin{Lem}\label{FP}
Suppose that \ $\bA \in \RR^{d \times d}_{(+)}$ \ is irreducible.
\renewcommand{\labelenumi}{{\rm(\roman{enumi})}}
\begin{enumerate}
 \item
  Then \ $s(\bA)$ \ is an eigenvalue of \ $\bA$, \ the algebraic and
   geometric multiplicities of \ $s(\bA)$ \ equal 1, and the real parts of the
   other eigenvalues of \ $\bA$ \ are less than \ $s(\bA)$.
 \item
  Corresponding to the eigenvalue \ $s(\bA)$ \ there exists a unique (right)
   eigenvector \ $\bu_\bA \in \RR^d_{++}$ \ of \ $\bA$ \  such that the sum of
   its coordinates is 1.
  The vector \ $\bu_\bA \in \RR^d_{++}$ \ is the unique (right)
   eigenvector of \ $\ee^\bA$ \ (called the right Perron vector of \ $\ee^\bA$)
   corresponding to the eigenvalue \ $r(\ee^\bA) = \ee^{s(\bA)}$ \ of \ $\ee^\bA$
   \ such that the sum of its coordinates is 1.
 \item
  There exists a unique left eigenvector \ $\bv_\bA \in \RR^d_{++}$ \ of \ $\bA$
   \ corresponding to the eigenvalue \ $s(\bA)$ \ with
   \ $\bu_\bA^\top \bv_\bA = 1$, \ and
   \[
     \ee^{-s(\bA)t} \ee^{t\bA} \to \bPi_\bA
     := \bu_\bA \bv_\bA^\top \in \RR^{d \times d}_{++} \qquad
     \text{as \ $t \to \infty$.}
   \]
  The vector \ $\bv_\bA \in \RR^d_{++}$ \ is the unique (left)
   eigenvector of \ $\ee^\bA$ \ (called the left Perron vector of \ $\ee^\bA$)
   corresponding to the eigenvalue \ $r(\ee^\bA) = \ee^{s(\bA)}$ \ of \ $\ee^\bA$
   \ such that \ $\bu_\bA^\top \bv_\bA = 1$.
 \item
  There exist \ $c_\bA, \kappa_\bA \in \RR_{++}$ \ such that
   \[
     \| \ee^{-s(\bA)t} \ee^{t\bA} - \bPi_\bA \| \leq c_\bA \ee^{- \kappa_\bA t} \qquad
     \text{for all \ $t \in \RR_+$.}
   \]
  Consequently,
   \ $\| \ee^{t\bA} \| \leq \bigr( c_\bA + \| \bPi_\bA \| \bigr) \ee^{-s(\bA)t}$,
   \ $t \in \RR_+$.
 \item
  Moreover,
   \[
     \frac{1}{t} \int_0^t \ee^{-s(\bA)u} \ee^{u\bA} \, \dd u \to \bPi_\bA \qquad
     \text{as \ $t \to \infty$.}
   \]
\end{enumerate}
\end{Lem}

\noindent
\textbf{Proof.}
The proof of (i) is based on the Frobenius--Perron theorem for \ $\ee^\bA$
 \ (see, e.g., Horn and Johnson \cite[Theorems 8.2.11 and 8.5.1]{HorJoh}).
Recall that \ $\sigma(\ee^\bA) = \ee^{\sigma(\bA)}$ \ and
 \ $r(\ee^\bA) = \ee^{s(\bA)}$, \ see \eqref{rs}.
By Lemma \ref{FP1}, \ $\ee^\bA \in \RR^{d \times d}_{++}$, \ hence, by the
 Frobenius--Perron theorem, \ $r(\ee^\bA) \in \RR_{++}$ \ is an eigenvalue of
 \ $\ee^\bA$, \ the algebraic and geometric multiplicities of \ $r(\ee^\bA)$
 \ equal 1, and the absolute values of the other eigenvalues of \ $\ee^\bA$
 \ are less than \ $r(\ee^\bA)$.
\ Then \ $s(\bA) = \log[r(\ee^\bA)]$ \ is an eigenvalue of \ $\bA$, \ the
 algebraic and geometric multiplicities of \ $s(\bA)$ \ equal 1, and the real
 parts of the other eigenvalues of \ $\bA$ \ are less than \ $s(\bA)$.
\ Indeed, since \ $\sigma(\ee^\bA) = \ee^{\sigma(\bA)}$, \ the algebraic
 multiplicity  of \ $s(\bA)$ \ (as an eigenvalue of \ $\bA$) \ coincides with
 the algebraic multiplicity of \ $r(\ee^\bA)$ \ (as an eigenvalue of
 \ $\ee^\bA$), \ yielding that the algebraic multiplicity of \ $s(\bA)$ \ is
 \ $1$; \ and using that the geometric multiplicity of \ $s(\bA)$ \ is less
 than or equal to its algebraic multiplicity, we obtain that the geometric
 multiplicity of  \ $s(\bA)$ \ equals \ $1$, \ too.
Further, if \ $\lambda$ \ is an eigenvalue of \ $\bA$ \ not equal to
 \ $s(\bA)$, \ then \ $\ee^\lambda$ \ is an eigenvalue of \ $\ee^\bA$ \ and
 \ $r(\ee^\bA) > |\ee^\lambda| = \ee^{\Re(\lambda)}$ \ yields that
 \ $s(\bA) = \log(r(\ee^\bA)) > \Re(\lambda)$, \ as it was stated, hence the
 proof of (i) is complete.

The Frobenius--Perron theorem also implies (ii) and the unique existence of
 \ $\bv_\bA$ \ in (iii), and the convergence in (iii) along the sequence of
 the positive integers.
The aim of the following discussion is to show that the convergence in (iii)
 holds also along the positive real numbers.
By Dunford and Schwartz \cite[Theorem VII.1.8]{DunSch},
 \[
   \ee^{t\bA}
   = \sum_{\lambda\in\sigma(\bA)}
      \sum_{i=0}^{\nu(\lambda)-1}
       \frac{(\bA - \lambda \bI_d)^i}{i!} t^i \ee^{\lambda t}
       E_\bA(\lambda) , \qquad t \in \RR_+ ,
 \]
 where \ $\nu(\lambda)$ \ denotes the index of \ $\lambda$ \ given in Dunford
 and Schwartz \cite[Definition VII.1.2]{DunSch}, and the projections
 \ $E_\bA(\lambda)$, \ $\lambda \in \sigma(\bA)$, \ are defined by
 \ $E_\bA(\lambda) := e_{\bA, \lambda}(\bA)$, \ where
 \ $e_{\bA, \lambda} : \CC \to \CC$ \ is an analytic function in some open set
 containing \ $\sigma(\bA)$ \ such that \ $e_{\bA, \lambda}(\mu) = 1$ \ if
 \ $\mu$ \ is in some neighborhood of \ $\lambda$ \ and
 \ $e_{\bA, \lambda}(\mu) = 0$ \ if \ $\mu$ \ is in some neighborhood of any
 point of \ $\sigma(\bA) \setminus \{\lambda\}$, \ see the definition of the
 function \ $e_{\bA, \lambda}$ \ in Dunford and Schwartz
 \cite[before Theorem VII.1.6]{DunSch}.
Here \ $e_{\bA, \lambda}(\bA) = P_{\bA,\lambda}(\bA)$, \ where
 \ $P_{\bA,\lambda} : \CC \to \CC$ \ is a polynomial with complex coefficients
 such that \ $P_{\bA,\lambda}^{(m)}(\mu) = e_{\bA, \lambda}^{(m)}(\mu)$ \ for all
 \ $\mu \in \sigma(\bA)$ \ and for all \ $m \in \{0, 1, \ldots, \nu(\mu)-1\}$,
 \ see the definition of the matrix \ $e_{\bA, \lambda}(\bA)$ \ in Dunford and
 Schwartz \cite[before Theorem VII.1.5]{DunSch}.
Consequently,
 \begin{equation}\label{DS}
  \ee^{t\bA}
  = \sum_{\lambda\in\sigma(\bA)}
     \sum_{i=0}^{\nu(\lambda)-1}
      \frac{(\bA - \lambda \bI_d)^i}{i!} t^i \ee^{\lambda t}
      P_{\bA,\lambda}(\bA) , \qquad t \in \RR_+ ,
 \end{equation}
 where \ $P_{\bA,\lambda}$ \ is a polynomial with complex coefficients such that
 \ $P_{\bA,\lambda}(\lambda) = 1$, \ $P_{\bA,\lambda}^{(m)}(\lambda) = 0$ \ for all
 \ $m \in \{1, \ldots, \nu(\lambda)-1\}$, \ and
 \ $P_{\bA,\lambda}^{(m)}(\tlambda) = 0$ \ for all \ $\tlambda \in \sigma(\bA)$
 \ with \ $\tlambda \ne \lambda$ \ and for all
 \ $m \in \{0, 1, \ldots, \nu(\tlambda)-1\}$.

Note that the index and the algebraic multiplicity of an eigenvalue
 \ $\lambda \in \sigma(\bA)$ \ coincide.
Indeed, the index \ $\nu(\lambda)$ \ is the smallest non-negative integer
 \ $\nu \in \ZZ_+$ \ such that \ $\cR_\lambda^{(\nu+1)} = \cR_\lambda^{(\nu)}$,
 \ where
 \[
   \cR_\lambda^{(\nu)}
   := \left\{ \bx \in \CC^d : (\bA - \lambda \bI_d)^\nu \bx = \bzero \right\} .
 \]
By a change of basis of \ $\RR^d$, \ we have \ $\bA = \bS \bJ \bS^{-1}$
 \ with some invertible matrix \ $\bS \in \CC^{d \times d}$, \ where \ $\bJ$
 \ denotes the Jordan normal form of \ $\bA$ \ consisting of Jordan blocks
 \ $\bJ_\mu \in \CC^{m(\mu) \times m(\mu)}$, \ $\mu \in \sigma(\bA)$, \ where
 \ $m(\mu)$ \ is the algebraic multiplicity of \ $\mu$.
\ One can easily verify that
 \ $\cR_\lambda^{(\nu)} = \bS \widetilde\cR_{\lambda}^{(\nu)}$ \ for all
 \ $\nu \in \ZZ_+$, \ where
 \[
   \widetilde\cR_\lambda^{(\nu)}
   := \left\{ \bx \in\CC^d : (\bJ - \lambda \bI_d)^\nu \bx = \bzero \right\} .
 \]
For each \ $\mu \in \sigma(\bA)$ \ with \ $\mu \ne \lambda$, \ we have
 \[
   \left\{ \bx \in \CC^{m(\mu)}
           : (\bJ_\mu - \lambda \bI_{m(\mu)})^k \bx = \bzero \right\}
   = \{ \bzero \}
 \]
 for all \ $k \in \NN$, \ since \ $(\bJ_\mu - \lambda \bI_{m(\mu)})^k$ \ is
 invertible.
Moreover, the dimension of the subspace
 \[
   \left\{ \bx \in \CC^{m(\lambda)}
           : (\bJ_\lambda - \lambda \bI_{m(\lambda)})^k \bx = \bzero \right\}
 \]
 equals \ $k$ \ for all \ $k \in \{1, \ldots, m(\lambda)\}$, \ hence
 \ $m(\lambda)$ \ is the smallest non-negative integer \ $\nu \in \ZZ_+$
 \ such that \ $\widetilde\cR_\lambda^{(\nu+1)} = \widetilde\cR_\lambda^{(\nu)}$.
\ Since \ $\cR_\lambda^{(\nu)} = \bS \widetilde\cR_{\lambda}^{(\nu)}$ \ for all
 \ $\nu \in \ZZ_+$, \ we conclude \ $\nu(\lambda) = m(\lambda)$.

By (i), the algebraic and geometric multiplicities of the eigenvalue
 \ $s(\bA)$ \ equal \ $1$, \ and hence we have \ $\nu(s(\bA)) = 1$.
\ Then \ $\ee^{-s(\bA)t} \ee^{t\bA} \to P_{\bA,s(\bA)}(\bA)$ \ as \ $t \to \infty$,
 \ since, for each \ $\lambda\in\sigma(\bA)$ \ with \ $\lambda \ne s(\bA)$
 \ we have \ $\Re(\lambda) < s(\bA)$, \ and then
 \ $|t^i \ee^{-s(\bA)t} \ee^{\lambda t}| =  t^i \ee^{-(s(\bA)-\Re(\lambda))t} \to 0$ \ as
 \ $t \to \infty$ \ for all \ $i\in\{0,1,\ldots,\nu(\lambda)-1\}$.
\ Applying the Frobenius--Perron theorem for the matrix \ $\ee^\bA$,
 \ we obtain
 \ $\ee^{-s(\bA)n} \ee^{n\bA} = r(\ee^\bA)^{-n} (\ee^\bA)^n \to \bPi_\bA
    = \bu_{\bA} \bv_{\bA}^\top$
 \ as \ $n \to \infty$ \ (along the sequence of the positive integers), thus
 \ $P_{\bA,s(\bA)}(\bA) = \bPi_\bA = \bu_\bA\bv_\bA^\top$, \ and we conclude (iii).

The statement (iv) is trivial for \ $d = 1$.
\ For \ $d \geq 2$, \ by formula \eqref{DS}, we obtain
 \begin{align*}
  \| \ee^{-s(\bA) t}\ee^{t\bA} - \bPi_\bA \|
  \leq \sum_{\lambda\in\sigma(\bA)\setminus\{s(\bA)\}}
        \sum_{i=0}^{\nu(\lambda)-1}
         \frac{\|\bA - \lambda \bI_d\|^i}{i!} t^i
         \| P_{\bA,\lambda}(\bA) \| \ee^{-(s(\bA) - \Re(\lambda))t}
 \end{align*}
 for all \ $t \in \RR_+$, \ hence we conclude (iv) with
 \begin{gather*}
  \kappa_\bA
  = \frac{1}{2}
    \left( s(\bA)
           - \max_{\lambda\in\sigma(\bA)\setminus\{s(\bA)\}} \Re(\lambda) \right)
  \in \RR_{++} , \\
  c_\bA = \sum_{\lambda\in\sigma(\bA)\setminus\{s(\bA)\}}
          \sum_{i=0}^{\nu(\lambda)-1}
           \frac{\|\bA - \lambda \bI_d\|^i}{i!}
           \| P_{\bA,\lambda}(\bA) \|
           \sup_{t\in\RR_+} (t^i \ee^{- \kappa_\bA t})
       \in \RR_{++} .
 \end{gather*}
Indeed, \ $c_\bA = 0$ \ would lead us to a contradiction, since then, by
 \ $\sup_{t\in\RR_+} (t^i \ee^{- \kappa_\bA t}) > 0$, \ for all
 \ $\lambda \in \sigma(\bA) \setminus \{s(\bA)\}$ \ and
 \ $i \in \{0, \ldots, \nu(\lambda)-1\}$, \ we get
 \ $\|\bA - \lambda \bI_d\|^i \, \| P_{\bA,\lambda}(\bA) \| = 0$, \ which implies
 \ $\| \ee^{-s(\bA) t}\ee^{t\bA} - \bPi_\bA \| = 0$ \ for all \ $t \in \RR_+$.
\ With a special choice of \ $t = 0$, \ this gives \ $\bPi_\bA = \bI_d$,
 \ which is a contradiction, since \ $\bPi_\bA \in \RR_{++}^{d\times d}$.
\ Further,
 \[
   \|\ee^{t\bA}\|
   = \ee^{s(\bA)t} \|\ee^{-s(\bA)t} \ee^{t\bA}\|
   \leq \ee^{s(\bA)t}
        \bigl( \|\ee^{-s(\bA)t} \ee^{t\bA} - \bPi_\bA\| + \|\bPi_\bA\| \bigr)
   \leq \bigl( c_\bA + \|\bPi_\bA\| \bigr) \ee^{s(\bA)t}
 \]
 for all \ $t \in \RR_+$.
\ Finally,
 \begin{align*}
  \left\| \frac{1}{t} \int_0^t \ee^{-s(\bA)u} \ee^{u\bA} \, \dd u
          - \bPi_\bA \right\|
  &= \left\| \frac{1}{t}
            \int_0^t (\ee^{-s(\bA)u} \ee^{u\bA} - \bPi_\bA) \, \dd u \right\|
  \leq \frac{1}{t} \int_0^t \|\ee^{-s(\bA)u} \ee^{u\bA} - \bPi_\bA\| \, \dd u \\
  &\leq \frac{1}{t} \int_0^t c_\bA \ee^{- \kappa_\bA u} \, \dd u
  = \frac{c_\bA}{\kappa_\bA \, t} (1 - \ee^{- \kappa_\bA t})
  \leq \frac{c_\bA}{\kappa_\bA \, t}
  \to 0
 \end{align*}
 as \ $t \to \infty$, \ which implies (v).
\proofend

\begin{Rem}
The Internet Seminar \cite{IS} contains some of the statements of Lemma
 \ref{FP}.
\proofend
\end{Rem}

\section{On moments of multi-type CBI processes}
\label{section_moments_CBI}

\begin{Pro}\label{mean_asymptotics}
 Let \ $(\bX_t)_{t\in\RR_+}$ \ be a multi-type CBI process with parameters
 \ $(d, \bc, \Bbeta, \bB, \nu, \bmu)$ \ such that \ $\EE(\|\bX_0\|) < \infty$
 \ and the moment condition \eqref{moment_condition_1} holds.
\ Suppose that \ $(\bX_t)_{t\in\RR_+}$ \ is irreducible.
Then the following assertions hold:
 \renewcommand{\labelenumi}{{\rm(\roman{enumi})}}
 \begin{enumerate}
  \item
   if \ $s(\tbB) < 0$, \ then
    \ $\lim_{t\to\infty} \EE(\bX_t) = - \tbB^{-1} \tBbeta$;
  \item
   if \ $s(\tbB) = 0$, \  then
    \ $\lim_{t\to\infty} t^{-1} \EE(\bX_t) = \bPi \tBbeta$;
  \item
   if \ $s(\tbB) > 0$, \ then
    \ $\lim_{t\to\infty} \ee^{-s(\tbB)t} \EE(\bX_t)
       = \bPi \EE(\bX_0) + \frac{1}{s(\tbB)} \bPi \tBbeta$,
  \end{enumerate}
 where \ $\bPi := \bPi_{\tbB} \in \RR_{++}^{d \times d}$ \ is defined in
 \textup{(iii)} of Lemma \ref{FP}.
\end{Pro}

\noindent
\textbf{Proof.}
If \ $s(\tbB) < 0$, \ then, by \eqref{EXbX},
 \[
   \lim_{t\to\infty} \EE(\bX_t)
   = \left( \int_0^\infty \ee^{u\tbB} \, \dd u \right) \tBbeta ,
 \]
 since the decomposition
 \begin{equation}\label{deco}
  \ee^{t\tbB}
  = \ee^{s(\tbB)t} \left( \ee^{-s(\tbB)t} \ee^{t\tbB} - \bPi \right)
    + \ee^{s(\tbB)t} \, \bPi ,
  \qquad t \in \RR_+ ,
 \end{equation}
 implies
 \[
   \|\ee^{t\tbB}\|
   \leq c \ee^{-(\kappa-s(\tbB))t} + \ee^{s(\tbB)t} \, \|\bPi\|
   \to 0 \qquad \text{as \ $t \to \infty$}
 \]
 and
 \begin{align*}
  \left\| \int_t^\infty \ee^{u\tbB} \, \dd u \right\|
  &\leq c \int_t^\infty \ee^{-(\kappa-s(\tbB))u} \, \dd u
        + \left( \int_t^\infty \ee^{s(\tbB)u} \, \dd u \right) \|\bPi\| \\
  &\leq \frac{c}{\kappa-s(\tbB)} \ee^{-(\kappa-s(\tbB))t}
        + \frac{1}{-s(\tbB)} \ee^{s(\tbB)t} \, \|\bPi\|
   \to 0 \qquad \text{as \ $t \to \infty$,}
 \end{align*}
 where \ $c := c_{\tbB} \in \RR_{++}$ \ and
 \ $\kappa := \kappa_{\tbB} \in \RR_{++}$ \ are given in (iv) of Lemma \ref{FP}.

Further,
 \[
   \tbB \int_0^t \ee^{u\tbB} \, \dd u
   = \ee^{t\tbB} - \bI_d
   \to - \bI_d \qquad \text{as \ $t \to \infty$,}
 \]
 hence \ $\int_0^\infty \ee^{u\tbB} \, \dd u = - \tbB^{-1}$, \ which yields (i).

If \ $s(\tbB) = 0$, \ then, again by the decomposition \eqref{deco},
 \ $\|\ee^{u\tbB}\| \leq c \ee^{-\kappa t} + \|\bPi\| \leq c + \|\bPi\|$ \ for
 all \ $t \in \RR_+$, \ thus \eqref{EXbX} and (v) of Lemma \ref{FP} imply
 (ii).

If \ $s(\tbB) > 0$, \ then the statement will follow
 from (iii) of Lemma \ref{FP} and
 \begin{equation}\label{(iii)}
  \lim_{t\to\infty} \ee^{-s(\tbB)t} \int_0^t \ee^{u\tbB} \, \dd u
  = \frac{1}{s(\tbB)} \bPi .
 \end{equation}
By the decomposition \eqref{deco},
 \begin{align*}
  &\left\| \ee^{-s(\tbB)t} \int_0^t \ee^{u\tbB} \, \dd u
          - \frac{1}{s(\tbB)} \bPi \right\| \\
  &\leq \left\| \ee^{-s(\tbB)t}
                \int_0^t
                 \ee^{s(\tbB)u}
                 \left( \ee^{-s(\tbB)u} \ee^{u\tbB} - \bPi \right)
                 \dd u \right\|
     + \left\| \ee^{-s(\tbB)t} \int_0^t \ee^{s(\tbB)u} \bPi \, \dd u
               - \frac{1}{s(\tbB)} \bPi \right\| \\
  &\leq c \ee^{-s(\tbB)t} \int_0^t \ee^{(s(\tbB)-\kappa)u} \, \dd u
        + \frac{\ee^{-s(\tbB)t}}{s(\tbB)} \|\bPi\|
   \leq \frac{c(\ee^{-\kappa t} - \ee^{-s(\tbB)t})}{s(\tbB) - \kappa}
        + \frac{\ee^{-s(\tbB)t}}{s(\tbB)} \|\bPi\|
   \to 0
 \end{align*}
 as \ $t \to \infty$, \ thus we obtain \eqref{(iii)}, and hence, by
 \eqref{EXbX}, (iii) as well.
\proofend

From Theorems 4.3 and 4.5 in Barczy et al. \cite{BarLiPap3}, we derive the
 following moment estimations for \ $(\bX_t)_{t\in\RR_+}$ \ and
 \ $(\bM_n)_{n\in\NN}$.

\begin{Lem}\label{moment_estimations_X_critical}
Let \ $(\bX_t)_{t\in\RR_+}$ \ be a multi-type CBI process with parameters
 \ $(d, \bc, \Bbeta, \bB, \nu, \bmu)$ \ such that \ $\EE(\|\bX_0\|^q) < \infty$
 \ and
 \begin{equation}\label{moment_condition_m}
  \int_{U_d} \|\bz\|^q \bbone_{\{\|\bz\|\geq1\}} \, \nu(\dd \bz) < \infty ,
  \qquad
  \int_{U_d} \|\bz\|^q \bbone_{\{\|\bz\|\geq1\}} \, \mu_i(\dd \bz) < \infty ,
  \quad
  i \in \{1, \ldots, d\}
 \end{equation}
 with some \ $q \in \NN$.
\ Suppose that \ $(\bX_t)_{t\in\RR_+}$ \ is irreducible and critical.
Then
 \begin{equation}\label{moment_ic}
  \sup_{t\in\RR_+} \frac{\EE(\|\bX_t\|^q)}{(1 + t)^q} < \infty .
 \end{equation}
In particular, \ $\EE(\|\bX_t\|^q) = \OO(t^q)$ \ as \ $t \to \infty$ \ in the
 sense that \ $\limsup_{t\to\infty} t^{-q} \EE(\|\bX_t\|^q) < \infty$.
\end{Lem}

\noindent
\textbf{Proof.}
By Theorem 4.3 in Barczy et al. \cite{BarLiPap3}, we have
 \ $\EE(\|\bX_t\|^q) < \infty$, \ $t \in \RR_+$, \ and
 \begin{align*}
  \EE(X_{t,j}^k)
  &\leq \tc(t)^k \EE(\|\bX_0\|^k)
        +  k \|\tBbeta\| \, \tc(t)^k
          \int_0^t \EE(\|\bX_s\|^{k-1}) \, \dd s \\
  &\quad
        + k (k - 1) \tc(t)^k
          \sum_{i=1}^d c_i \int_0^t \EE(\|\bX_s\|^{k-1}) \, \dd s \\
  &\quad
        + \tc(t)^k
          \sum_{\ell=0}^{k-2}
           \binom{k}{\ell}
           \biggl[ \sum_{i=1}^d
                   \int_0^t \EE(\|\bX_s\|^{\ell+1}) \, \dd s
                   \int_{U_d}
                    \|\bz\|^{k-\ell} \, \mu_i(\dd \bz) \\
  &\phantom{\leq + \tc(t)^k \sum_{\ell=0}^{k-2} \binom{k}{\ell} \biggl[ \:}
                   + \int_0^t \EE(\|\bX_s\|^\ell) \, \dd s
                     \int_{U_d}
                      \|\bz\|^{k-\ell} \, \nu(\dd \bz) \biggr]
 \end{align*}
 for all \ $k \in \{1, \ldots, q\}$, \ $t \in \RR_+$, \ and
 \ $j \in \{1, \ldots, d\}$, \ where \ $\bX_t:=(X_{t,j})_{j\in\{1,\ldots,d\}}$,
 \ and \ $\tc(t) := \sup_{u\in[0,t]} \|\ee^{u\tbB}\| \leq c + \|\bPi\|$ \ for all
 \ $t \in \RR_+$ \ due to \eqref{c+Pi}, since \ $(\bX_t)_{t\in\RR_+}$ \ is
 irreducible and critical.
We will show \eqref{moment_ic} by induction with respect to
 \ $k \in \{1, \ldots, q\}$.
\ If \ $k = 1$, \ then the above estimate implies
 \begin{align*}
  \EE(X_{t,j})
  \leq \tc(t) \EE(\|\bX_0\|) + \|\tBbeta\| \tc(t) t
  \leq (c + \|\bPi\|) (\EE(\|\bX_0\|) + \|\tBbeta\|t) , \qquad
  t \in \RR_+ .
 \end{align*}
This yields \ $\sup_{t\in\RR_+} (1 + t)^{-1} \EE(\|\bX_t\|) < \infty$.
\ Further, for all \ $k \in \{2, \ldots, q\}$, \ $q\geq 2$, \ and
 \ $\ell \in \{0, 1, \ldots, k-2\}$, \ we have
 \begin{align*}
  \int_{U_d} \|\bz\|^{k-\ell} \, \nu(\dd \bz)
  &= \int_{U_d} \|\bz\|^{k-\ell} \bbone_{\{\|z\|<1\}} \, \nu(\dd \bz)
     + \int_{U_d} \|\bz\|^{k-\ell} \bbone_{\{\|z\|\geq1\}} \, \nu(\dd \bz) \\
  &\leq \int_{U_d} \|\bz\|^2 \bbone_{\{\|z\|<1\}} \, \nu(\dd \bz)
        + \int_{U_d} \|\bz\|^k \bbone_{\{\|z\|\geq1\}} \, \nu(\dd \bz)
   < \infty
 \end{align*}
 by \eqref{moment_condition_2_m} (for this we use \ $q\geq 2$) \ and
 \eqref{moment_condition_m}.
The finiteness of \ $\int_{U_d} \|\bz\|^{k-\ell} \, \mu_i(\dd \bz)$,
 \ $i \in \{1, \ldots, d\}$, \ follows in the same way for all
 \ $k \in \{2, \ldots, q\}$, \ $q \geq 2$, \ and
 \ $\ell \in \{0, 1, \ldots, k-2\}$.
\ Hence, by the power means inequality, the above estimate yields
 \ $\sup_{t\in\RR_+} (1 + t)^{-k} \EE(\|\bX_t\|^k) < \infty$.
\proofend

\begin{Lem}\label{moment_estimations_1_2}
Let \ $(\bX_t)_{t\in\RR_+}$ \ be a multi-type CBI process with parameters
 \ $(d, \bc, \Bbeta, \bB, \nu, \bmu)$ \ such that \ $\EE(\|\bX_0\|^{2q}) < \infty$
 \ and the moment conditions \eqref{moment_condition_m} hold with some
 \ $q \in \NN$.
\ Suppose that \ $(\bX_t)_{t\in\RR_+}$ \ is irreducible and critical.
Then, for the martingale differences
 \ $\bM_n = \bX_n - \EE(\bX_n \mid \bX_{n-1})$, \ $n \in \NN$, \ we have
 \ $\EE(\|\bM_n\|^{2q}) = \OO(n^q)$ \ as \ $n \to \infty$, \ i.e.,
 \ $\sup_{n\in \NN} n^{-q} \EE(\|\bM_n\|^{2q})<\infty$.
\end{Lem}

\noindent
\textbf{Proof.}
Applying Theorem 4.5 in Barczy et al. \cite{BarLiPap3}, we obtain
 \begin{align*}
  \EE(M_{n,j}^{2q} \mid \bX_{n-1} = \bx)
  &= \EE\left[ (X_{n,j} - \EE(X_{n,j} \mid \bX_{n-1}))^{2q}
               \mid \bX_{n-1} = \bx \right] \\
  &= \EE\left[ (X_{1,j} - \EE(X_{1,j} \mid \bX_0))^{2q}
               \mid \bX_0 = \bx \right]
   = P_{1,2q,j}(\bx)
 \end{align*}
 for all \ $n \in \NN$, \ $\bx \in \RR_+^d$, \ and \ $j \in \{1, \ldots, d\}$,
 \ where \ $\bM_n:=(M_{n,j})_{j\in\{1,\ldots,d\}}$, \ and hence
 \[
   \EE(M_{n,j}^{2q} \mid \bX_{n-1}) = P_{1,2q,j}(\bX_{n-1}) ,
 \]
 where \ $P_{1,2q,j} : \RR^d \to \RR$ \ is a polynomial having degree at most
 \ $q$.
\ Using Lemma \ref{moment_estimations_X_critical}, this yields
 \[
   \EE(M_{n,j}^{2q})
   = \EE\left[ P_{1,2q,j}(\bX_{n-1}) \right]
   = \OO(n^q)
 \]
 for all \ $j \in \{1, \ldots, d\}$.
By the power means inequality,
 \[
    \left(\frac{1}{d}\sum_{k=1}^d a_k^2 \right)^{\frac{1}{2}}
       \leq \left(\frac{1}{d}\sum_{k=1}^d a_k^{2q} \right)^{\frac{1}{2q}},
     \qquad a_1,\ldots,a_d\in\RR,
 \]
 and hence \ $(a_1^2+\cdots+a_d^2)^q\leq d^{q-1}(a_1^{2q}+\cdots+a_d^{2q})$, \ $a_1,\ldots,a_d\in\RR$.
\ Then we conclude
 \[
   \EE(\|\bM_n\|^{2q})
     =\EE\left(\left(\sum_{k=1}^d M_{n,j}^2\right)^q\right)
     \leq d^{q-1} \sum_{k=1}^d \EE(M_{n,j}^{2q})
     = \OO(n^{q}).
 \]
\proofend

\section{Convergence of random step processes}
\label{section_conv_step_drocesses}

Next we recall a result about convergence of random step processes towards a
 diffusion process, see Isp\'any and Pap \cite[Corollary 2.2]{IspPap2}.

\begin{Thm}\label{Conv2DiffThm}
Let \ $\bgamma : \RR_+ \times \RR^d \to \RR^{d \times r}$ \ be a continuous
 function.
Assume that uniqueness in the sense of probability law holds for the SDE
 \begin{equation}\label{SDE}
  \dd \, \bcU_t
  = \gamma (t, \bcU_t) \, \dd \bcW_t ,
  \qquad t \in \RR_+,
 \end{equation}
 with initial value \ $\bcU_0 = \bu_0$ \ for all \ $\bu_0 \in \RR^d$, \ where
 \ $(\bcW_t)_{t \in \RR_+}$ \ is an $r$-dimensional Brownian motion.
Let \ $(\bcU_t)_{t \in \RR_+}$ \ be a solution of the SDE \eqref{SDE} with
 initial value \ $\bcU_0 = \bzero$.

For each \ $n \in \NN$, \ let \ $(\bU^{(n)}_k)_{k \in \ZZ_+}$ \ be a sequence of
 \ $d$-dimensional martingale differences with respect to a filtration
 \ $(\cF^{(n)}_k)_{k \in \ZZ_+}$ \ (i.e.,
 \ $\EE(\bU_k^{(n)} \mid \cF_{k-1}^{(n)}) = \bU_{k-1}^{(n)}$, \  $k \in \NN$).
\ Let
 \[
   \bcU^{(n)}_t := \sum_{k=0}^{\nt} \bU^{(n)}_k  \, ,
   \qquad t \in \RR_+, \quad n \in \NN .
 \]
Suppose \ $\EE \big( \|\bU^{(n)}_k\|^2 \big) < \infty$ \ for all
 \ $n, k \in \NN$, \ and \ $\bU_0^{(n)} \distr \bzero$ \ as \ $n \to \infty$.
\ Suppose that for each \ $T > 0$,
 \begin{enumerate}
  \item[\textup{(i)}]
   $\sup\limits_{t\in[0,T]}
     \left\| \sum\limits_{k=1}^{\nt}
              \var\Bigl( \bU^{(n)}_k \mid \cF^{(n)}_{k-1} \Bigr)
             - \int_0^t
                \bgamma(s,\bcU^{(n)}_s) \bgamma(s,\bcU^{(n)}_s)^\top
                \dd s \right\|
         \stoch 0$,\\
  \item[\textup{(ii)}]
   $\sum\limits_{k=1}^{\lfloor nT \rfloor}
     \EE \big( \|\bU^{(n)}_k\|^2 \bbone_{\{\|\bU^{(n)}_k\| > \theta\}}
                 \bmid \cF^{(n)}_{k-1} \big)
    \stoch 0$
   \ for all \ $\theta > 0$,
 \end{enumerate}
 where \ $\stoch$ \ denotes convergence in probability.
Then \ $\bcU^{(n)} \distr \bcU$ \ as \ $n \to \infty$.
\end{Thm}

\section{A version of the continuous mapping theorem}
\label{continuous_mapping_theorem}

For functions \ $f$ \ and \ $f_n$, \ $n \in \NN$, \ in \ $\DD(\RR_+, \RR^d)$,
 \ we write \ $f_n \lu f$ \ if \ $(f_n)_{n \in \NN}$ \ converges to \ $f$
 \ locally uniformly, i.e., if \ $\sup_{t \in [0,T]} \|f_n(t) - f(t)\| \to 0$
 \ as \ $n \to \infty$ \ for all \ $T \in \RR_{++}$.
\ For measurable mappings \ $\Phi : \DD(\RR_+, \RR^d) \to \DD(\RR_+, \RR^q)$
 \ and \ $\Phi_n : \DD(\RR_+, \RR^d) \to \DD(\RR_+, \RR^q)$, \ $n \in \NN$,
 \ we will denote by \ $C_{\Phi, (\Phi_n)_{n \in \NN}}$ \ the set of all functions
 \ $f \in \CC(\RR_+, \RR^d)$ \ for which \ $\Phi_n(f_n) \to \Phi(f)$ \ whenever
 \ $f_n \lu f$ \ with \ $f_n \in \DD(\RR_+, \RR^d)$, \ $n \in \NN$.

\begin{Lem}\label{Conv2Funct}
Let \ $(\bcU_t)_{t \in \RR_+}$ \ and \ $(\bcU^{(n)}_t)_{t \in \RR_+}$, \ $n \in \NN$,
 \ be \ $\RR^d$-valued stochastic processes with c\`adl\`ag paths such that
 \ $\bcU^{(n)} \distr \bcU$ \ as \ $n\to\infty$.
\ Let \ $\Phi : \DD(\RR_+, \RR^d) \to \DD(\RR_+, \RR^q)$ \ and
 \ $\Phi_n : \DD(\RR_+, \RR^d) \to \DD(\RR_+, \RR^q)$, \ $n \in \NN$, \ be
 measurable mappings such that there exists
 \ $C \subset C_{\Phi,(\Phi_n)_{n\in\NN}}$ \ with \ $C \in \cD_\infty(\RR_+, \RR^d)$
 \ and \ $\PP(\bcU \in C) = 1$.
\ Then \ $\Phi_n(\bcU^{(n)}) \distr \Phi(\bcU)$ \ as \ $n\to\infty$.
\end{Lem}

Lemma \ref{Conv2Funct} can be considered as a consequence of Theorem 3.27 in
 Kallenberg \cite{Kal}, and we note that a proof of this lemma can also be
 found in Isp\'any and Pap \cite[Lemma 3.1]{IspPap2}.

\end{document}